\documentclass{amsart}

\usepackage{amssymb}
\usepackage{amsmath}
\usepackage{amsthm}
\usepackage{amsfonts}

\usepackage{a4wide}

\newtheorem{proposition}{Proposition}
\newtheorem{theorem}{Theorem}
\newtheorem{corollary}{Corollary}

\theoremstyle{definition}

\newtheorem{example}{Example}

\newtheorem{remark}{Remark}
\newtheorem{conjecture}{Conjecture}

\begin{document}

%\linenumbers

\title[Elevated rank distribution]
{BCF-groups with elevated rank distribution}

\author{Daniel C. Mayer}
\address{Naglergasse 53\\8010 Graz\\Austria}
\email{algebraic.number.theory@algebra.at}
\urladdr{http://www.algebra.at}

\thanks{Research supported by the Austrian Science Fund (FWF): projects J0497-PHY, P26008-N25, and by EUREA}

\subjclass[2010]{Primary 20D15, 20E18, 20E22, 20F05, 20F12, 20F14;
secondary 11R37, 11R32, 11R11, 11R20, 11R29, 11Y40
}

\keywords{pro-\(3\) groups, finite \(3\)-groups, BCF-groups,
generator rank, relation rank, Schur \(\sigma\)-groups, balanced presentation, extremal root path principle,
low index normal subgroups, kernels of Artin transfers,
abelian quotient invariants of first and second order, elevated rank distribution,
\(p\)-group generation algorithm, descendant trees, antitony principle;
three-stage Hilbert \(3\)-class field towers, maximal unramified pro-\(3\) extensions,
unramified cyclic cubic extensions, unramified nonic extensions,
imaginary quadratic fields, non-elementary bicyclic \(3\)-class groups, Galois action,
punctured capitulation type B.18, minimal discriminants}

\date{Thursday, 07 October 2021}

%--------------------------------------------------------------------------------

\begin{abstract}
Infinitely many large Schur \(\sigma\)-groups \(G\) with logarithmic order \(\mathrm{lo}(G)=19+e\),
non-elementary bicyclic commutator quotient \(G/G^\prime\simeq C_{3^e}\times C_3\), \(e\ge 2\),
elevated rank distribution \(\varrho(G)=(3,3,3;3)\),
punctured transfer kernel type \(\varkappa(G)\sim (144;4)\)
and soluble length \(\mathrm{sl}(G)=3\) are constructed.
Up to \(e\le 4\), they are realized as \(3\)-class field tower groups \(\mathrm{Gal}(\mathrm{F}_3^\infty(K)/K)\)
of imaginary quadratic number fields \(K=\mathbb{Q}(\sqrt{d})\), \(d<0\).
Their metabelianizations \(M=G/G^{\prime\prime}\) are BCF-groups
with \(\mathrm{lo}(M)=8+e\) and
bicyclic third lower central factor \(\gamma_3(M)/\gamma_4(M)\simeq C_3\times C_3\).
\end{abstract}

\maketitle

%\newpage
%--------------------------------------------------------------------------------

\section{Introduction}
\label{s:Intro}

\noindent
Let \(G\) be a pro-\(3\) group or finite \(3\)-group
with bicyclic commutator quotient
\(G/G^\prime\simeq C_{3^e}\times C_3\)
having one non-elementary factor with exponent \(e\ge 2\).
Then \(G\) possesses four maximal self-conjugate subgroups
\(H_1,\ldots,H_3;H_4\),
and by the \textit{rank distribution} of \(G\) we understand the quartet
\begin{equation}
\label{eqn:Rho}
\varrho(G):=\left\lbrack\mathrm{rank}_3(H_1/H_1^\prime),\ldots,\mathrm{rank}_3(H_3/H_3^\prime);\mathrm{rank}_3(H_4/H_4^\prime)\right\rbrack.
\end{equation}
Let \((\gamma_j(G))_{j\ge 1}\) be the lower central series of \(G\).
When the factors \(\gamma_j(G)/\gamma_{j-1}(G)\simeq C_3\) are all cyclic,
for \(j\ge 2\),
then \(G\) is called a \textit{CF-group}, according to Ascione et al.
\cite{AHL1977}.
CF means cyclic factors.
Otherwise, at least the factor
\(\gamma_3(G)/\gamma_4(G)\simeq C_3\times C_3\) is bicyclic,
and \(G\) is called a \textit{BCF-group}, according to Nebelung
\cite{Ne1989}.
BCF means bicyclic or cyclic factors.
Recall that, since \(\gamma_2(G)=\langle s_2,\gamma_3(G)\rangle\),
the factor \(\gamma_2(G)/\gamma_3(G)\simeq C_3\) is always cyclic,
generated by the main commutator \(s_2=\lbrack y,x\rbrack\)
of the two-generated group \(G=\langle x,y\rangle\).
For a BCF-group \(G\), we have
\(\gamma_3(G)=\langle s_3,t_3,\gamma_4(G)\rangle\)
with higher non-trivial commutators
\(s_3=\lbrack s_2,x\rbrack\) and \(t_3=\lbrack s_2,y\rbrack\).

%--------------------------------------------------------------------------------

In
\cite[\S\ 2]{Ma2021},
we introduced the concept of \textit{punctured transfer kernel types}
\begin{equation}
\label{eqn:Kappa}
\varkappa(G):=\left\lbrack\ker(T_1),\ldots,\ker(T_3);\ker(T_4)\right\rbrack
\end{equation}
for \(3\)-groups \(G=\langle x,y\rangle\) with
\(G/G^\prime\simeq C_{3^e}\times C_3\), \(e\ge 2\).
Here, \(T_i:\,G/G^\prime\to H_i/H_i^\prime\)
denotes the Artin transfer homomorphism from \(G\) to \(H_i\).
It turned out that at least three kernels are two-dimensional,
equal to the complete \(3\)-elementary subgroup
\(\langle x^{e-1},y,G^\prime\rangle/G^\prime\) of \(G/G^\prime\),
when \(G\) is a CF-group.
Consequently, metabelian CF-groups can be realized arithmetically
only by second \(3\)-class groups \(\mathrm{Gal}(\mathrm{F}_3^2(K)/K)\)
of real quadratic fields \(K=\mathbb{Q}(\sqrt{d})\), \(d>0\),
but not for imaginary quadratic fields with \(d<0\),
where all kernels must be one-dimensional.

%--------------------------------------------------------------------------------

In
\cite[\S\S\ 5 and 7]{Ma2021},
we investigated how BCF-groups \(G\) with \textit{moderate} rank distribution
\(\varrho(G)\in\lbrace (2,2,2;3),(2,2,3;3)\rbrace\)
are populated by second \(3\)-class groups \(\mathrm{Gal}(\mathrm{F}_3^2(K)/K)\)
and \(3\)-class field tower groups \(\mathrm{Gal}(\mathrm{F}_3^\infty(K)/K)\)
of imaginary quadratic fields \(K=\mathbb{Q}(\sqrt{d})\), \(d<0\),
with non-elementary bicyclic \(3\)-class groups
\(\mathrm{Cl}_3(K)\simeq C_{3^e}\times C_3\), \(e\ge 2\).

In the present article we continue this research enterprise
for BCF-groups \(G\) with \textit{elevated} rank distribution
\(\varrho(G)=(3,3,3;3)\)
and punctured transfer kernel type \(\mathrm{B}.18\), \(\varkappa(G)\sim (144;4)\).
Their exo-genetic propagation has been clarified in
\cite[Thm. 17]{Ma2021}.

%\newpage
%--------------------------------------------------------------------------------

\section{Arithmetical realization}
\label{s:NumberTheory}

\noindent
It is of the greatest importance to emphasize that
the assumptions concerning the punctured transfer kernel type
\(\varkappa(G)\),
and the \textit{logarithmic abelian quotient invariants} of first order
\begin{equation}
\label{eqn:Alpha1}
\alpha_1(G):=\left\lbrack H_1/H_1^\prime,\ldots,H_3/H_3^\prime;H_4/H_4^\prime\right\rbrack
\end{equation}
and of second order
\begin{equation}
\label{eqn:Alpha2}
\alpha_2(G):=\left(G/G^\prime;\lbrack H_i/H_i^\prime;(H_{i,j}/H_{i,j}^\prime)_{(H_i:H_{i,j})=3}\rbrack_{1\le i\le 4}\right)
\end{equation}
of the Schur \(\sigma\)-groups in the following six main theorems
are perfectly tailored for applications
in algebraic number theory and class field theory.
According to the Artin reciprocity law
\cite{Ar1927,Ar1929},
these invariants can be interpreted
for an arbitrary algebraic number field \(K\)
as the punctured capitulation type
\(\varkappa(K):=\left\lbrack\ker(\tau_1),\ldots,\ker(\tau_3);\ker(\tau_4)\right\rbrack\)
of the extension homomorphisms
\(\tau_i:\,\mathrm{Cl}_3(K)\to\mathrm{Cl}_3(L_i)\),
\(\mathfrak{a}\mathcal{P}_K\mapsto(\mathfrak{a}\mathcal{O}_{L_i})\mathcal{P}_{L_i}\),
of \(3\)-classes from \(K\) to the four unramified cyclic cubic extensions \(L_i\),
the logarithmic abelian type invariants
\(\alpha_1(K):=\left\lbrack\mathrm{Cl}_3(L_1),\ldots,\mathrm{Cl}_3(L_3);\mathrm{Cl}_3(L_4)\right\rbrack\)
of the \(3\)-class groups of the fields \(L_i\),
and the logarithmic abelian type invariants
\[\alpha_2(K):=\left(\mathrm{Cl}_3(K);\lbrack\mathrm{Cl}_3(L_i);(\mathrm{Cl}_3(L_{i,j}))_{\lbrack L_{i,j}:L_i\rbrack=3}\rbrack_{1\le i\le 4}\right)\]
of all unramified (but not necessarily abelian) \(3\)-extensions of degree at most nine of \(K\).
For details see
\cite{Ma2020a}.
In this article,
we investigate applications to the simplest algebraic number fields,
namely imaginary quadratic fields \(K=\mathbb{Q}(\sqrt{d})\)
with negative fundamental discriminants \(d<0\).

%\newpage
%--------------------------------------------------------------------------------

\section{Main theorems}
\label{s:MainTheorems}

\noindent
The six \textit{main theorems} are the crucial achievements of the present article.
They show that Schur \(\sigma\)-groups
\cite{KoVe1975,Ag1998,BBH2017}
with \textit{elevated} rank distribution
also become \textit{periodic} for sufficiently large exponents \(e\ge 9\),
similarly as Schur \(\sigma\)-groups with \textit{moderate} rank distribution
for \(e\ge 5\), according to
\cite{Ma2021}.

\subsection{Schur \(\sigma\)-groups \(G\)}
\label{ss:SchurSigmaGroups}

\noindent
In all theorems, the symbols \(e^+:=e+1\), \(e^-:=e-1\) are used for abbreviation.
Isomorphism classes of groups are identified in accordance with
\cite{BEO2005,GNO2006,MAGMA2021}.

\begin{theorem}
\label{thm:Elevated24}
A total of \(54\) Schur \(\sigma\)-groups \(G\) with
commutator quotient \(G/G^\prime\simeq (3^e,3)\),
punctured transfer kernel type \(\mathrm{B}.18\), \(\varkappa(G)\sim (144;4)\),
elevated rank distribution  \(\varrho(G)=(3,3,3;3)\),
first abelian quotient invariants \(\alpha_1(G)\sim\lbrack e^+21,e11,e11;e^-21\rbrack\),
second abelian quotient invariants
\begin{equation}
\label{eqn:Elevated24AQI2}
\begin{aligned}
\alpha_2(G)\sim (e1;&\lbrack e^+21;e2111,(e^+211)^3,(e^+2)^9 \rbrack, \\
&\lbrack e11;e2111,(e^+21)^3,(e21)^9 \rbrack, \\
&\lbrack e11;e2111,(e^+111)^3,(e^+2)^9 \rbrack; \\
&\lbrack e^-21;e2111,(e31)^3,(e21)^8,e^-22 \rbrack )
\end{aligned}
\end{equation}
and (minimal) logarithmic order \(\mathrm{lo}(G)=19+e\) is given for each \(e\ge 9\) by the term
\begin{equation}
\label{eqn:Elevated24SchurSigma}
G=W_a\lbrack-\#1;1\rbrack^{e-9}-\#1;p-\#1;q-\#1;1 \text{ with } p\in\lbrace 2,3\rbrace \text{ and } q\in\lbrace 1,2,3\rbrace,
\end{equation}
where \(9\) distinct periodic roots with
\(1\le a\le 9\), \(\tilde{a}=1\) for \(a\le 3\), \(\tilde{a}=2\) for \(a\ge 4\),
are denoted by
\begin{equation}
\label{eqn:Elevated24PeriodicRoots}
W_a:=\langle 2187,3\rangle-\#3;2-\#\mathbf{4;24}-\#3;14-\#4;a-\#2;\tilde{a}-\#2;1.
\end{equation}
\end{theorem}

%--------------------------------------------------------------------------------

\begin{theorem}
\label{thm:Elevated26}
A total of \(162\) Schur \(\sigma\)-groups \(G\) with
commutator quotient \(G/G^\prime\simeq (3^e,3)\),
punctured transfer kernel type \(\mathrm{B}.18\), \(\varkappa(G)\sim (144;4)\),
elevated rank distribution  \(\varrho(G)=(3,3,3;3)\),
first abelian quotient invariants \(\alpha_1(G)\sim\lbrack e^+21,e11,e11;e^-21\rbrack\),
second abelian quotient invariants
\begin{equation}
\label{eqn:Elevated26AQI2}
\begin{aligned}
\alpha_2(G)\sim (e1;&\lbrack e^+21;e2111,(e^+1111)^3,(e^+2)^9 \rbrack, \\
&\lbrack e11;e2111,(e^+21)^3,(e21)^9 \rbrack, \\
&\lbrack e11;e2111,(e^+21)^3,(e^+2)^9 \rbrack; \\
&\lbrack e^-21;e2111,(e31)^3,(e21)^8,e^-22 \rbrack )
\end{aligned}
\end{equation}
and (minimal) logarithmic order \(\mathrm{lo}(G)=19+e\) is given for each \(e\ge 9\) by the term
\begin{equation}
\label{eqn:Elevated26SchurSigma}
G=W_{a,b}\lbrack-\#1;1\rbrack^{e-9}-\#1;p-\#1;q-\#1;1 \text{ with } p\in\lbrace 2,3\rbrace \text{ and } 1\le q\le N,
\end{equation}
where \(45\) distinct periodic roots with
\(1\le a\le 27\), \(\tilde{a}=1\), \(1\le b\le 3\), \(N=1\) for \(a\in\lbrace 3,4,8,\) \(12,13,17,21,22,26\rbrace\)
and \(\tilde{a}=2\), \(b=1\), \(N=3\) otherwise,
are denoted by
\begin{equation}
\label{eqn:Elevated26PeriodicRoots}
W_{a,b}:=\langle 2187,3\rangle-\#3;2-\#\mathbf{4;26}-\#3;14-\#4;a-\#2;\tilde{a}-\#2;b.
\end{equation}
\end{theorem}

%--------------------------------------------------------------------------------

\begin{theorem}
\label{thm:Elevated28or30}
A total of \(324\) Schur \(\sigma\)-groups \(G\) with
commutator quotient \(G/G^\prime\simeq (3^e,3)\),
punctured transfer kernel type \(\mathrm{B}.18\), \(\varkappa(G)\sim (144;4)\),
elevated rank distribution  \(\varrho(G)=(3,3,3;3)\),
first abelian quotient invariants \(\alpha_1(G)\sim\lbrack e^+21,e11,e11;e^-21\rbrack\),
second abelian quotient invariants
\begin{equation}
\label{eqn:Elevated28or30AQI2}
\begin{aligned}
\alpha_2(G)\sim (e1;&\lbrack e^+21;e2111,(e^+211)^3,(e^+2)^9 \rbrack, \\
&\lbrack e11;e2111,(e^+21)^3,(e^+2)^9 \rbrack, \\
&\lbrack e11;e2111,(e^+111)^3,(e^+2)^9 \rbrack; \\
&\lbrack e^-21;e2111,(e31)^3,(e21)^8,e^-22 \rbrack )
\end{aligned}
\end{equation}
and (minimal) logarithmic order \(\mathrm{lo}(G)=19+e\) is given for each \(e\ge 9\) by the term
\begin{equation}
\label{eqn:Elevated28or30SchurSigma}
G=W_{\ell,k,a}\lbrack-\#1;1\rbrack^{e-9}-\#1;p-\#1;q-\#1;r \text{ with } p\in\lbrace 2,3\rbrace \text{ and } q,r\in\lbrace 1,2,3\rbrace,
\end{equation}
where \(18\) distinct periodic roots with
\((\ell,k)\in\lbrace \mathbf{(28,5),(30,2)}\rbrace\) and \(1\le a\le 9\)
are denoted by
\begin{equation}
\label{eqn:Elevated28or30PeriodicRoots}
W_{\ell,k,a}:=\langle 2187,3\rangle-\#3;2-\mathbf{\#4;\ell-\#3;k}-\#4;a-\#2;1-\#2;1.
\end{equation}
\end{theorem}

%--------------------------------------------------------------------------------

\begin{theorem}
\label{thm:Elevated31}
A total of \(162\) Schur \(\sigma\)-groups \(G\) with
commutator quotient \(G/G^\prime\simeq (3^e,3)\),
punctured transfer kernel type \(\mathrm{B}.18\), \(\varkappa(G)\sim (144;4)\),
elevated rank distribution  \(\varrho(G)=(3,3,3;3)\),
first abelian quotient invariants \(\alpha_1(G)\sim\lbrack e^+21,e11,e11;e^-21\rbrack\),
second abelian quotient invariants
\begin{equation}
\label{eqn:Elevated31AQI2}
\begin{aligned}
\alpha_2(G)\sim (e1;&\lbrack e^+21;e2111,(e^+211)^3,(e^+2)^9 \rbrack, \\
&\lbrack e11;e2111,(e^+21)^3,(e^+2)^9 \rbrack, \\
&\lbrack e11;e2111,(e^+21)^3,(e21)^9 \rbrack; \\
&\lbrack e^-21;e2111,(e211)^3,(e21)^8,e^-22 \rbrack )
\end{aligned}
\end{equation}
and (minimal) logarithmic order \(\mathrm{lo}(G)=19+e\) is given for each \(e\ge 9\) by the term
\begin{equation}
\label{eqn:Elevated31SchurSigma}
G=W_{a,b}\lbrack-\#1;1\rbrack^{e-9}-\#1;p-\#1;q-\#1;1 \text{ with } p\in\lbrace 2,3\rbrace \text{ and } 1\le q\le N,
\end{equation}
where \(45\) distinct periodic roots with
\(1\le a\le 27\), \(\tilde{a}=1\) for \(a\le 9\), \(\tilde{a}=2\) for \(a\ge 10\),
\(1\le b\le 3\), \(N=1\) for \(a\in\lbrace 1,5,9,11,15,16,21,22,26\rbrace\) and \(b=1\), \(N=3\) otherwise,
are denoted by
\begin{equation}
\label{eqn:Elevated31PeriodicRoots}
W_{a,b}:=\langle 2187,3\rangle-\#3;2-\#\mathbf{4;31}-\#3;29-\#4;a-\#2;\tilde{a}-\#2;b.
\end{equation}
\end{theorem}

%--------------------------------------------------------------------------------

\begin{theorem}
\label{thm:Elevated33}
A total of \(162\) Schur \(\sigma\)-groups \(G\) with
commutator quotient \(G/G^\prime\simeq (3^e,3)\),
punctured transfer kernel type \(\mathrm{B}.18\), \(\varkappa(G)\sim (144;4)\),
elevated rank distribution  \(\varrho(G)=(3,3,3;3)\),
first abelian quotient invariants \(\alpha_1(G)\sim\lbrack e^+21,e11,e11;e^-21\rbrack\),
second abelian quotient invariants
\begin{equation}
\label{eqn:Elevated33AQI2}
\begin{aligned}
\alpha_2(G)\sim (e1;&\lbrack e^+21;e2111,(e^+1111)^3,(e^+2)^9 \rbrack, \\
&\lbrack e11;e2111,(e^+21)^3,(e^+2)^9 \rbrack, \\
&\lbrack e11;e2111,(e^+21)^3,(e^+2)^9 \rbrack; \\
&\lbrack e^-21;e2111,(e31)^3,(e21)^8,e^-22 \rbrack )
\end{aligned}
\end{equation}
and (minimal) logarithmic order \(\mathrm{lo}(G)=19+e\) is given for each \(e\ge 9\) by the term
\begin{equation}
\label{eqn:Elevated33SchurSigma}
G=W_{a,b}\lbrack-\#1;1\rbrack^{e-9}-\#1;p-\#1;q-\#1;1 \text{ with } p\in\lbrace 2,3\rbrace \text{ and } q\in\lbrace 1,2,3\rbrace,
\end{equation}
where \(27\) distinct periodic roots with
\(1\le a\le 9\) and \(1\le b\le 3\)
are denoted by
\begin{equation}
\label{eqn:Elevated33PeriodicRoots}
W_{a,b}:=\langle 2187,3\rangle-\#3;2-\#\mathbf{4;33}-\#3;32-\#4;a-\#2;1-\#2;b.
\end{equation}
\end{theorem}

%--------------------------------------------------------------------------------

\begin{theorem}
\label{thm:Elevated37}
A total of \(162\) Schur \(\sigma\)-groups \(G\) with
commutator quotient \(G/G^\prime\simeq (3^e,3)\),
punctured transfer kernel type \(\mathrm{B}.18\), \(\varkappa(G)\sim (144;4)\),
elevated rank distribution  \(\varrho(G)=(3,3,3;3)\),
first abelian quotient invariants \(\alpha_1(G)\sim\lbrack e^+21,e11,e11;e^-21\rbrack\),
second abelian quotient invariants
\begin{equation}
\label{eqn:Elevated37AQI2}
\begin{aligned}
\alpha_2(G)\sim (e1;&\lbrack e^+21;e2111,(e^+211)^3,(e^+2)^9 \rbrack, \\
&\lbrack e11;e2111,(e^+21)^3,(e^+2)^9 \rbrack, \\
&\lbrack e11;e2111,(e^+21)^3,(e^+2)^9 \rbrack; \\
&\lbrack e^-21;e2111,(e211)^3,(e21)^8,e^-22 \rbrack )
\end{aligned}
\end{equation}
and (minimal) logarithmic order \(\mathrm{lo}(G)=19+e\) is given for each \(e\ge 9\) by the term
\begin{equation}
\label{eqn:Elevated37SchurSigma}
G=W_{a,b}\lbrack-\#1;1\rbrack^{e-9}-\#1;p-\#1;q-\#1;1 \text{ with } p\in\lbrace 2,3\rbrace \text{ and } q\in\lbrace 1,2,3\rbrace,
\end{equation}
where \(27\) periodic roots with
\(1\le a\le 9\), \(\tilde{a}=1\) for \(a\in\lbrace 2,6,7\rbrace\), \(\tilde{a}=2\) otherwise, and \(1\le b\le 3\)
are
\begin{equation}
\label{eqn:Elevated37PeriodicRoots}
W_{a,b}:=\langle 2187,3\rangle-\#3;2-\#\mathbf{4;37}-\#3;32-\#4;a-\#2;\tilde{a}-\#2;b.
\end{equation}
\end{theorem}

%--------------------------------------------------------------------------------

\begin{remark}
\label{rmk:Elevated}
The \textit{periodic twig} \(-\#1;p-\#1;q-\#1;r\) of the terms
for the Schur \(\sigma\)-groups \(G\) in the main theorems
contains \(6\) \textit{terminal leaves} on average.
However, for Theorem
\ref{thm:Elevated28or30}
there are \(18\),
and for Theorems
\ref{thm:Elevated26}
and
\ref{thm:Elevated31}
there are partially only \(2\).

For each \(e\ge 9\),
all main theorems together
yield \(1026\) Schur \(\sigma\)-groups \(G\)
with \(\mathrm{lo}(G)=19+e\),
which are descendants of \(171\) distinct periodic roots \(W\)
with fixed logarithmic order \(\mathrm{lo}(W)=25\).
\end{remark}

%\newpage
%--------------------------------------------------------------------------------

\noindent
Exemplarily we give a succinct proof for the last main theorem,
namely Theorem
\ref{thm:Elevated37}.

\begin{proof}
(Proof of Theorem
\ref{thm:Elevated37}.)
For a fixed step size \(s\ge 1\),
we denote by \(N\) the number of all immediate descendants of a \(3\)-group, and
by \(C\) the number of capable immediate descendants with positive nuclear rank \(\nu\ge 1\).
Generally, let \(X:=\langle 2187,3\rangle-\#3;2-\#4;37-\#3;32\).
This is a non-metabelian \(3\)-group of type \((729,3)\).
We consider a chain of exo-genetic propagations:
\begin{itemize}
\item
\(X\) has \(N=C=27\) for \(s=\nu=4\)
but only the first \(9\) descendants are of type \((2187,3)\).
\item
Each \(X-\#4;a\) with \(1\le a\le 9\)
has \(N=C=6\) for \(s=\nu=2\)
but only the first, resp. second, descendant,
indicated by \(\tilde{a}\in\lbrace 1,2\rbrace\), is of type \((6561,3)\).
\item
Each \(X-\#4;a-\#2,\tilde{a}\) with \(1\le a\le 9\)
has \(N=C=9\) for \(s=\nu=2\)
but only the first \(3\) descendants are of type \((19683,3)\).
\item
Each \(W_{a,b}:=X-\#4;a-\#2,\tilde{a}-\#2;b\) with \(1\le a\le 9\) and \(1\le b\le 3\)
has \(6\) Schur \(\sigma\)-descendants
\(W_{a,b}\lbrack-\#1;1\rbrack^{e-9}-\#1;p-\#1;q-\#1;1\) with \(p\in\lbrace 2,3\rbrace\) and \(q\in\lbrace 1,2,3\rbrace\),
for each \(e\ge 9\).
\end{itemize}
Together this census yields \(9\cdot 3\cdot 6=162\) Schur \(\sigma\)-groups, for each \(e\ge 9\).
\end{proof}

%--------------------------------------------------------------------------------

The following supplementary theorem
provides a warranty for the fact
that the information in the six main theorems is exhaustive and complete.

\begin{theorem}
\label{thm:Exhaustion}
(\textbf{Exhaustion Theorem.}) \\
Let \(G\) be a Schur \(\sigma\)-group with
non-elementary bicyclic commutator quotient \(G/G^\prime\simeq C_{3^e}\times C_3\), \(e\ge 9\),
punctured transfer kernel type \(\mathrm{B}.18\), \(\varkappa(G)\sim (144;4)\),
elevated rank distribution \(\varrho(G)=(3,3,3;3)\), and
first abelian quotient invariants \(\alpha_1(G)\sim\lbrack (e+1)21,e11,e11;(e-1)21\rbrack\).
Then
\begin{itemize}
\item
if \(G\) has logarithmic order \(\mathrm{lo}(G)=19+e\),
then \(G\) is of one of the shapes in the six main Theorems
\ref{thm:Elevated24} --
\ref{thm:Elevated37}
(inclusively the shape of \(\alpha_2(G)\))
and has soluble length \(\mathrm{sl}(G)=3\);
\item
if \(G\) is not of one of the shapes in the six main Theorems
\ref{thm:Elevated24} --
\ref{thm:Elevated37},
then \(G\) has logarithmic order \(\mathrm{lo}(G)>19+e\)
and different second abelian quotient invariants \(\alpha_2(G)\).
\end{itemize}
\end{theorem}

%\newpage
%--------------------------------------------------------------------------------

\subsection{Second derived quotients \(G/G^{\prime\prime}\)}
\label{ss:Metabelianizations}

\noindent
Periodicity of metabelianizations with elevated rank distribution
sets in earlier for \(e\ge 5\) already.

\begin{corollary}
\label{cor:Elevated}
The metabelianization \(M=G/G^{\prime\prime}\)
of a Schur \(\sigma\)-group \(G\) with
commutator quotient \(G/G^\prime\simeq (3^e,3)\), \(e\ge 5\),
punctured transfer kernel type \(\mathrm{B}.18\), \(\varkappa(G)\sim (144;4)\),
logarithmic abelian quotient invariants of first order
\(\alpha_1(G)\sim\lbrack (e+1)21,e11,e11;(e-1)21\rbrack\),
and logarithmic order \(\mathrm{lo}(G)=19+e\)
is given by one of the two candidates
\begin{equation}
\label{eqn:ElevatedMetabelianization}
M\simeq\langle 2187,3\rangle-\#3;2-\#2;93\lbrack-\#1;1\rbrack^{e-5}-\#1;i \text{ with } i\in\lbrace 2,3\rbrace.
\end{equation}
Their logarithmic order is \(\mathrm{lo}(M)=8+e\),
i.e. the second derived subgroup \(G^{\prime\prime}\)
is of constant logarithmic order \(\mathrm{lo}(G^{\prime\prime})=11\),
in fact, it is abelian of constant type \(G^{\prime\prime}\simeq (332111)\).
A parametrized polycyclic power commutator presentation
of the members \(\langle 2187,3\rangle-\#3;2-\#2;93\lbrack-\#1;1\rbrack^{e-5}\)
of the infinite chain is given for \(e\ge 6\) by
\begin{equation}
\label{eqn:Presentation}
\langle x,y\mid x^{3^e}=1, y^3=s_3s_4^2, s_2^3=s_4t_4^2, s_3^3=s_5, t_3^3=s_5^2, \lbrack x^3,y\rbrack=s_4t_4s_5^2, \lbrack x^3,s_2\rbrack=s_5, t_5=s_5\rangle
\end{equation}
in terms of the commutators
\(s_2=\lbrack y,x\rbrack\), 
\(s_3=\lbrack s_2,x\rbrack\), \(t_3=\lbrack s_2,y\rbrack\), 
\(s_4=\lbrack s_3,x\rbrack\), \(t_4=\lbrack t_3,y\rbrack\), 
\(s_5=\lbrack s_4,x\rbrack\), \(t_5=\lbrack t_4,y\rbrack\). 
\end{corollary}

\noindent
The justification of the periodicities in Theorems
\ref{thm:Elevated24} --
\ref{thm:Elevated37}
and Corollary
\ref{cor:Elevated}
will be developed in \S\
\ref{s:Periodicities}.

%\newpage
%--------------------------------------------------------------------------------

\section{Layout of the paper}
\label{s:Layout}

\noindent
Since the periodicity in the crucial Theorems
\ref{thm:Elevated24} --
\ref{thm:Elevated37}
sets in with exponent \(e=9\),
we devote \S\S\
\ref{s:21},
\ref{s:31}
and
\ref{s:41}
to the detailed discussion of the regular cases \(2\le e\le 4\).
We do not go into the details of the irregular intermediate cases
\(5\le e\le 8\),
which are clarified sufficiently by Figure
\ref{fig:SchurSigma59}.
In \S\
\ref{s:Periodicities}
we illuminate the long and winding road to the actual verification of the periodicity
of Schur \(\sigma\)-groups \(G\) with elevated rank distribution \(\varrho(G)=(3,3,3;3)\)
and commutator quotient \(G/G^\prime\simeq (3^e,3)\),
which was expected by ourselves for \(e\ge 9\) in analogy to the periodicity
for \(e\ge 5\) in the case of moderate rank distribution
\cite{Ma2021}.
Arithmetical applications to \(3\)-class field tower groups
\(\mathrm{Gal}(\mathrm{F}_3^\infty(K)/K)\)
of imaginary quadratic fields \(K=\mathbb{Q}(\sqrt{d})\)
with fundamental discriminants \(d<0\)
and non-elementary \(3\)-class groups \(\mathrm{Cl}_3(K)\simeq (3^e,3)\)
are given in \S\S\
\ref{s:Imag21},
\ref{s:Imag31}
and
\ref{s:Imag41}
for \(2\le e\le 4\). In \S\S\
\ref{s:Imag51}
and
\ref{s:Imag61},
where arithmetical realizations of \(5\le e\le 6\) are just possible (with CPU-time a week)
it becomes clear that \(e=7\) (CPU-time several months) and \(e=8\) (CPU-time several years)
are outside of a reasonable and realistic arithmetical enterprise,
aggravated by internal Magma errors, due to huge absolute discriminants \(\lvert d\rvert\).
A conclusion concerning the general structure of the
logarithmic abelian quotient invariants \(\alpha_2\) of second order
is eventually drawn in \S\
\ref{s:General}.
In \S\S\
\ref{s:21},
\ref{s:31}
and
\ref{s:41},
we also consider
\(\mathrm{lo}(G)>19+e\).
The case \(e=2\) was also investigated in
\cite{Ma2020b}.

%\newpage
%--------------------------------------------------------------------------------

\section{Root path to Schur \(\sigma\)-groups}
\label{s:RootPath}

\noindent
In order to find \(\sigma\)-groups
\cite[Dfn. 3.1, p. 91]{Ma2018},
and in particular Schur \(\sigma\)-groups
\cite{KoVe1975,Ag1998,BBH2017},
\(G\) with commutator quotient \(G/G^\prime\simeq (3^e,3)\)
and punctured transfer kernel type \(\mathrm{B}.18\), \(\varkappa(G)\sim (144;4)\),
it is necessary to take into consideration the associated \textit{scaffold type} 
\(\mathrm{b}.31\), \(\varkappa\sim (044;4)\),
since the two-dimensional transfer kernel \(0\) of a parent
can shrink to the one-dimensional transfer kernel \(1\) for a descendant.
This is a consequence of the \textit{antitony principle} for the Artin pattern \((\varkappa,\alpha)\)
of parent descendant pairs.
The situation is similar to
\cite[\S\ 3.2.2 and Fig. 2, pp. 91--92]{Ma2018}
and
\cite[Fig. 1--2, pp. 24--25]{Ma2015b},
both for elementary \(G/G^\prime\simeq (3,3)\).
Now we have non-elementary \(G/G^\prime\).

\begin{proposition}
\label{prp:RootPath}
The root path of the \textbf{bifurcation} \(B:=\langle 2187,3\rangle-\#3;2\) \textbf{of infinite order},
\begin{equation}
\label{eqn:RootPath}
1\buildrel s=2 \over\longleftarrow\pi_p^3(B)=\langle 9,2\rangle\buildrel s=2 \over\longleftarrow\pi_p^2(B)=\langle 81,3\rangle
\buildrel s=3 \over\longleftarrow\pi_p(B)=\langle 2187,3\rangle\buildrel s=3 \over\longleftarrow B=\langle 2187,3\rangle-\#3;2,
\end{equation}
has step sizes \((2,2,3,3)\) and
contains two vertices with scaffold type \(\mathrm{b}.31\),
\(\varkappa\sim (044;4)\),
which give rise to Schur \(\sigma\)-groups \(G\) with type \(\mathrm{B}.18\),
\(\varkappa(G)\sim (144;4)\),
and to their metabelianizations \(G/G^{\prime\prime}\).
\end{proposition}

\begin{proof}
There are only three groups \(G\) with \(G/G^\prime\simeq (9,3)\), i.e. \(e=2\),
and order \(\#G=81\), namely the non-abelian groups
\(G\simeq\langle 81,3\rangle\) with \(\varkappa(G)\sim (000;0)\), \(\mathrm{a}.1\),
\(G\simeq\langle 81,4\rangle\) with \(\varkappa(G)\sim (444;4)\), \(\mathrm{A}.20\), and
\(G\simeq\langle 81,6\rangle\) with \(\varkappa(G)\sim (111;1)\), \(\mathrm{A}.1\).
According to the antitony principle for the Artin pattern \((\varkappa,\alpha)\),
the latter two groups are discouraged as predecessors of descendants with
\(\varkappa\sim (044;4)\) or \(\varkappa\sim (144;4)\).
Moreover, they are not \(\sigma\)-groups.
The unique remaining group \(G=\langle 81,3\rangle\) has the root path
\(G\buildrel s=2 \over\longrightarrow\pi_p(G)=\langle 9,2\rangle=C_3\times C_3\buildrel s=2 \over\longrightarrow\pi_p^2(G)=\langle 1,1\rangle=1\).
In order to stay at \(e=2\),
the descendant \(D=\langle 729,10\rangle\) with scaffold type \(\mathrm{b}.31\), \(\varkappa(D)\sim (044;4)\),
must be selected.
The unique immediate \(\sigma\)-descendant \(F=\langle 6561,165\rangle\) of \(D\)
is already the fork between the desired Schur \(\sigma\)-group \(S\)
and its metabelianization \(S/S^{\prime\prime}\simeq F-\#2;85\).
See \S\
\ref{s:Periodicities},
 Figure
\ref{fig:SchurSigma}.
\end{proof}

%\newpage
%--------------------------------------------------------------------------------

\section{\(3\)-groups with commutator quotient \((9,3)\)}
\label{s:21}

\noindent
In Table
\ref{tbl:21},
we list the second AQI \(\alpha_2\)
of the \(30\) non-metabelian step size-\(4\) descendants
\(F-\#4;\ell\) with \(1\le\ell\le 30\)
of the metabelian fork \(F=\langle 729,10\rangle-\#2;2\).
The general structure of \(\alpha_2\) is
\begin{equation}
\label{eqn:2AQI21}
\alpha_2(G)=\lbrack 21;(\tau_0;22111,D_1),(211;22111,D_2),(211;22111,D_3);(211;22111,D_4)\rbrack,
\end{equation}
where each dodecuplet \(D_i\), \(1\le i\le 4\), consists of a triplet \(T_i^3\) and a nonet \(N_i^9\).
The metabelianization \(M=G/G^{\prime\prime}\) is given by the step size-\(2\) descendant
\(F-\#2;m\) with \(m\in\lbrace 82,83,84,85\rbrace\) of \(F\).
The smallest logarithmic order, soluble length, of a Schur \(\sigma\)-descendant \(S\) of \(G\)
is \(\mathrm{lo}(S)\), \(\mathrm{sl}(S)\).

%\newpage
%--------------------------------------------------------------------------------

\renewcommand{\arraystretch}{1.0}
\begin{table}[ht]
\caption{Invariants of \(G=\langle 729,10\rangle-\#2;2-\#4;\ell\) with \(1\le\ell\le 30\)}
\label{tbl:21}
\begin{center}
\begin{tabular}{|r|c||l|l||l|l||l|l||l|l||c||c|c|}
\hline
 \(\ell\) & \(\tau_0\) & \(T_1\)   & \(N_1\) & \(T_2\)  & \(N_2\) & \(T_3\)  & \(N_3\) & \(T_4\)  & \(N_4\)  & \(m\)  & \(\mathrm{lo}(S)\) & \(\mathrm{sl}(S)\) \\
\hline
    \(1\) & \(222\)    & \(22111\) & \(221\) & \(2211\) & \(221\) & \(2211\) & \(221\) & \(3111\) &  \(32\)  & \(82\) & \(\infty\) & \(\infty\) \\
\hline
    \(2\) & \(222\)    & \(22111\) & \(221\) & \(321\)  & \(32\)  & \(321\)  & \(32\)  & \(321\)  &  \(32\)  & \(83\) & \(21\) & \(3\) \\
    \(7\) & \(222\)    & \(22111\) & \(221\) & \(321\)  & \(32\)  & \(321\)  & \(32\)  & \(321\)  &  \(32\)  & \(82\) & \(21\) & \(3\) \\
\hline
    \(3\) & \(222\)    & \(3211\)  & \(221\) & \(321\)  & \(32\)  & \(321\)  & \(32\)  & \(3111\) &  \(32\)  & \(82\) & \(25\) & \(4\) \\
\hline
    \(4\) & \(222\)    & \(3211\)  & \(221\) & \(321\)  & \(32\)  & \(321\)  & \(32\)  & \(2211\) &  \(221\) & \(83\) & \(21\) & \(3\) \\
    \(6\) & \(222\)    & \(3211\)  & \(221\) & \(321\)  & \(32\)  & \(321\)  & \(32\)  & \(2211\) &  \(221\) & \(83\) & \(21\) & \(3\) \\
\hline
    \(5\) & \(222\)    & \(3211\)  & \(221\) & \(321\)  & \(32\)  & \(321\)  & \(32\)  & \(3111\) &  \(32\)  & \(82\) & \(25\) & \(3\) \\
\hline
    \(8\) & \(222\)    & \(22111\) & \(221\) & \(321\)  & \(221\) & \(321\)  & \(221\) & \(321\)  &  \(32\)  & \(83\) & \(24\) & \(4\) \\
\hline
    \(9\) & \(222\)    & \(22111\) & \(221\) & \(3111\) & \(32\)  & \(3111\) & \(32\)  & \(3111\) &  \(32\)  & \(83\) & \(\infty\) & \(\infty\) \\
\hline
   \(10\) & \(222\)    & \(3211\)  & \(221\) & \(321\)  & \(32\)  & \(321\)  & \(221\) & \(3111\) &  \(32\)  & \(82\) & \(21\) & \(3\) \\
   \(12\) & \(222\)    & \(3211\)  & \(221\) & \(321\)  & \(32\)  & \(321\)  & \(221\) & \(3111\) &  \(32\)  & \(83\) & \(21\) & \(3\) \\
   \(13\) & \(222\)    & \(3211\)  & \(221\) & \(321\)  & \(32\)  & \(321\)  & \(221\) & \(3111\) &  \(32\)  & \(82\) & \(21\) & \(3\) \\
   \(15\) & \(222\)    & \(3211\)  & \(221\) & \(321\)  & \(32\)  & \(321\)  & \(221\) & \(3111\) &  \(32\)  & \(83\) & \(21\) & \(3\) \\
\hline
   \(11\) & \(222\)    & \(3211\)  & \(221\) & \(321\)  & \(32\)  & \(321\)  & \(32\)  & \(3111\) &  \(32\)  & \(83\) & \(21\) & \(3\) \\
   \(14\) & \(222\)    & \(3211\)  & \(221\) & \(321\)  & \(32\)  & \(321\)  & \(32\)  & \(3111\) &  \(32\)  & \(83\) & \(21\) & \(3\) \\
\hline
   \(16\) & \(321\)    & \(3211\)  & \(32\)  & \(321\)  & \(32\)  & \(321\)  & \(32\)  & \(3111\) &  \(32\)  & \(84\) & \(25\) & \(4\) \\
\hline
   \(17\) & \(321\)    & \(3211\)  & \(32\)  & \(321\)  & \(32\)  & \(321\)  & \(221\) & \(3111\) &  \(32\)  & \(85\) & \(21\) & \(3\) \\
   \(24\) & \(321\)    & \(3211\)  & \(32\)  & \(321\)  & \(32\)  & \(321\)  & \(221\) & \(3111\) &  \(32\)  & \(85\) & \(21\) & \(3\) \\
\hline
   \(18\) & \(321\)    & \(3211\)  & \(32\)  & \(321\)  & \(221\) & \(321\)  & \(221\) & \(3111\) &  \(32\)  & \(84\) & \(28\) & \(4\) \\
\hline
   \(19\) & \(321\)    & \(3211\)  & \(32\)  & \(321\)  & \(32\)  & \(321\)  & \(32\)  & \(3111\) &  \(32\)  & \(85\) & \(21\) & \(3\) \\
   \(23\) & \(321\)    & \(3211\)  & \(32\)  & \(321\)  & \(32\)  & \(321\)  & \(32\)  & \(3111\) &  \(32\)  & \(85\) & \(21\) & \(3\) \\
   \(25\) & \(321\)    & \(3211\)  & \(32\)  & \(321\)  & \(32\)  & \(321\)  & \(32\)  & \(3111\) &  \(32\)  & \(84\) & \(21\) & \(3\) \\
   \(27\) & \(321\)    & \(3211\)  & \(32\)  & \(321\)  & \(32\)  & \(321\)  & \(32\)  & \(3111\) &  \(32\)  & \(85\) & \(21\) & \(3\) \\
\hline
   \(20\) & \(321\)    & \(31111\) & \(32\)  & \(3111\) & \(32\)  & \(3111\) & \(32\)  & \(3111\) &  \(32\)  & \(84\) & \(\infty\) & \(\infty\) \\
\hline
   \(21\) & \(321\)    & \(31111\) & \(32\)  & \(321\)  & \(32\)  & \(321\)  & \(32\)  & \(321\)  &  \(221\) & \(85\) & \(21\) & \(3\) \\
   \(28\) & \(321\)    & \(31111\) & \(32\)  & \(321\)  & \(32\)  & \(321\)  & \(32\)  & \(321\)  &  \(221\) & \(84\) & \(21\) & \(3\) \\
\hline
   \(22\) & \(321\)    & \(3211\)  & \(32\)  & \(321\)  & \(32\)  & \(321\)  & \(32\)  & \(2211\) &  \(221\) & \(84\) & \(21\) & \(3\) \\
\hline
   \(26\) & \(321\)    & \(3211\)  & \(32\)  & \(321\)  & \(32\)  & \(321\)  & \(221\) & \(2211\) &  \(221\) & \(85\) & \(24\) & \(4\) \\
\hline
   \(29\) & \(321\)    & \(31111\) & \(32\)  & \(321\)  & \(32\)  & \(321\)  & \(32\)  & \(321\)  &  \(32\)  & \(85\) & \(21\) & \(3\) \\
\hline   
   \(30\) & \(321\)    & \(31111\) & \(32\)  & \(3111\) & \(32\)  & \(3111\) & \(32\)  & \(2211\) &  \(221\) & \(85\) & \(\infty\) & \(\infty\) \\
\hline
\end{tabular}
\end{center}
\end{table}

%--------------------------------------------------------------------------------

\begin{theorem}
\label{thm:21}
The Schur \(\sigma\)-groups \(S\) with commutator quotient \(S/S^\prime\simeq (9,3)\),
punctured transfer kernel type \(\mathrm{B}.18\), \(\varkappa(S)\sim (144;4)\),
and first AQI \(\alpha_1(S)\sim (\tau_0,211,211;211)\)
are descendants of \(30\) non-metabelian \(3\)-groups
\(G=\langle 729,10\rangle-\#2;2-\#4;\ell\) whose invariants are listed in Table
\ref{tbl:21}.
In the case of finite order \(\mathrm{lo}(S)<\infty\),
their invariants usually coincide with those of the predecessor \(G\).
For \(\mathrm{lo}(S)=21\) they have three stages, \(\mathrm{sl}(S)=3\),
for \(\mathrm{lo}(S)\in\lbrace 24,28\rbrace\) four stages, \(\mathrm{sl}(S)=4\),
and for \(\mathrm{lo}(S)=25\) they have \(3\le\mathrm{sl}(S)\le 4\).
Their metabelianization \(S/S^{\prime\prime}\simeq G/G^{\prime\prime}\)
is \(M=\langle 729,10\rangle-\#2;2-\#2;m\), where
\(m\in\lbrace 82,83\rbrace\), \(\tau_0=222\) for \(1\le\ell\le 15\), and
\(m\in\lbrace 84,85\rbrace\), \(\tau_0=321\) for \(16\le\ell\le 30\).
\end{theorem}

%--------------------------------------------------------------------------------

\noindent
The minimum \(\mathrm{lo}(S)=21\) occurs for \(20\) values \(\ell\),
\(24\) for \(2\), \(25\) for \(3\), \(28\) for \(1\), and \(\infty\) for \(4\).

%\newpage
%--------------------------------------------------------------------------------

\section{Imaginary quadratic fields \(K\) with \(\mathrm{Cl}_3(K)\simeq C_9\times C_3\)}
\label{s:Imag21}

\noindent
The \(875\) imaginary quadratic fields \(K=\mathbb{Q}(\sqrt{d})\)
with fundamental discriminants \(-1\,000\,000<d<0\) and
\(3\)-class group \(\mathrm{Cl}_3(K)\simeq C_9\times C_3\)
were computed by means of Magma
\cite{MAGMA2021}
in \(7\,782\) seconds of CPU time, that is more than two hours.
In Table
\ref{tbl:Type21},
the first nineteen cases with punctured capitulation type \(\mathrm{B}.18\),
\(\varkappa(K)\sim (144;4)\),
are listed.
The abelian quotient invariants \(\alpha_1(K)\) of first order
of only eleven of them are \textit{uni-polarized} and in the \textit{ground state}.
For details see
\cite{Ma2020b}.

%\newpage
%--------------------------------------------------------------------------------

\renewcommand{\arraystretch}{1.0}
\begin{table}[ht]
\caption{Nineteen fields \(K=\mathbb{Q}(\sqrt{d})\) with \(\mathrm{Cl}_3(K)\simeq C_9\times C_3\) and \(\varkappa(K)\sim (144;4)\)}
\label{tbl:Type21}
\begin{center}
\begin{tabular}{|r|r|c||c|c|}
\hline
 No.     & \(d\)         & factors          & \(\alpha_1(K)\)       & remark \\
\hline
  \(45\) &  \(-89\,923\) & prime            & \((222,211,211;321)\) & bi-polarized \\
  \(87\) & \(-150\,319\) & \(13,31,373\)    & \((321,211,211;211)\) & \\
 \(124\) & \(-194\,703\) & \(3,64\,901\)    & \((321,211,211;211)\) & \\
 \(161\) & \(-242\,255\) & \(5,13,3\,727\)  & \((222,211,211;321)\) & bi-polarized \\
 \(203\) & \(-294\,983\) & \(13,22\,691\)   & \((222,211,211;211)\) & \\
 \(304\) & \(-389\,371\) & \(401,971\)      & \((431,211,211;211)\) & first excited state \\
 \(305\) & \(-389\,435\) & \(5,71,1\,097\)  & \((222,211,211;211)\) & \\
 \(330\) & \(-409\,380\) & \(2,3,5,6\,823\) & \((222,211,211;321)\) & bi-polarized \\
 \(397\) & \(-481\,567\) & \(271,1\,777\)   & \((222,211,211;321)\) & bi-polarized \\
 \(413\) & \(-494\,771\) & \(61,8\,111\)    & \((321,211,211;211)\) & \\
 \(418\) & \(-497\,859\) & \(3,263,631\)    & \((321,211,211;321)\) & bi-polarized \\
 \(438\) & \(-518\,835\) & \(3,5,34\,589\)  & \((222,211,211;211)\) & \\
 \(470\) & \(-553\,807\) & \(433,1\,279\)   & \((222,211,211;211)\) & \\
 \(482\) & \(-566\,168\) & \(2,17,23,181\)  & \((321,211,211;211)\) & \\
 \(635\) & \(-761\,855\) & \(5,17,8\,963\)  & \((222,211,211;211)\) & \\
 \(637\) & \(-763\,972\) & \(2,11,97,179\)  & \((222,211,211;211)\) & \\
 \(661\) & \(-793\,992\) & \(2,3,33\,083\)  & \((321,211,211;321)\) & bi-polarized \\
 \(729\) & \(-857\,743\) & prime            & \((431,211,211;321)\) & highly bi-polarized \\
 \(743\) & \(-876\,948\) & \(2,3,73\,079\)  & \((222,211,211;211)\) & \\
\hline
\end{tabular}
\end{center}
\end{table}

%\newpage
%--------------------------------------------------------------------------------

\noindent
In Table
\ref{tbl:SecondAQI21},
we give the abelian quotient invariants \(\alpha_2(K)\) of second order
of the eleven fields in the uni-polarized ground state
contained in Table
\ref{tbl:Type21}.
The general structure of \(\alpha_2(K)\) is the following
\begin{equation}
\label{eqn:SecondAQI21}
\alpha_2(K)=\lbrack 21;(\tau_0;22111,D_1),(211;22111,D_2),(211;22111,D_3);(211;22111,D_4)\rbrack
\end{equation}
where \(\tau_0\in\lbrace 222,321\rbrace\), and
each dodecuplet \(D_i\), \(1\le i\le 4\), consists of a triplet and a nonet.

%\newpage
%--------------------------------------------------------------------------------

\renewcommand{\arraystretch}{1.1}
\begin{table}[ht]
\caption{Details for eleven fields \(K=\mathbb{Q}(\sqrt{d})\) in Table \ref{tbl:Type21}}
\label{tbl:SecondAQI21}
\begin{center}
\begin{tabular}{|r||c|c|c|c|c|c|}
\hline
 No.     & \(\tau_0\) & \(D_1\)              & \(D_2\)            & \(D_3\)             & \(D_4\)             & remark \\
\hline
  \(87\) & \(321\)    & \((41111)^3(32)^9\)  & \((321)^3(32)^9\)  & \((321)^3(32)^9\)   & \((321)^3(32)^9\)   & ref. 29 \\
 \(124\) & \(321\)    & \((3211)^3(32)^9\)   & \((321)^3(32)^9\)  & \((321)^3(32)^9\)   & \((3111)^3(32)^9\)  & ref. 16,19,23,25,27 \\
 \(203\) & \(222\)    & \((32111)^3(221)^9\) & \((3211)^3(32)^9\) & \((3211)^3(221)^9\) & \((2221)^3(221)^9\) & extreme \\
 \(305\) & \(222\)    & \((32111)^3(221)^9\) & \((3211)^3(32)^9\) & \((3211)^3(221)^9\) & \((2221)^3(221)^9\) & extreme \\
 \(413\) & \(321\)    & \((3211)^3(32)^9\)   & \((321)^3(32)^9\)  & \((321)^3(32)^9\)   & \((3111)^3(32)^9\)  & ref. 16,19,23,25,27 \\
 \(438\) & \(222\)    & \((3211)^3(221)^9\)  & \((321)^3(32)^9\)  & \((321)^3(32)^9\)   & \((2211)^3(221)^9\) & ref. 4,6 \\
 \(470\) & \(222\)    & \((3211)^3(221)^9\)  & \((321)^3(32)^9\)  & \((321)^3(221)^9\)  & \((3111)^3(32)^9\)  & ref 10,12,13,15 \\
 \(482\) & \(321\)    & \((32211)^3(32)^9\)  & \((3221)^3(32)^9\) & \((3221)^3(32)^9\)  & \((3221)^3(221)^9\) & extreme \\
 \(635\) & \(222\)    & \((3211)^3(221)^9\)  & \((321)^3(32)^9\)  & \((321)^3(32)^9\)   & \((3111)^3(32)^9\)  & ref 3,5,11,14 \\
 \(637\) & \(222\)    & \((3211)^3(221)^9\)  & \((321)^3(32)^9\)  & \((321)^3(221)^9\)  & \((3111)^3(32)^9\)  & ref 10,12,13,15 \\
 \(743\) & \(222\)    & \((3211)^3(221)^9\)  & \((321)^3(32)^9\)  & \((321)^3(221)^9\)  & \((3111)^3(32)^9\)  & ref 10,12,13,15 \\
\hline
\end{tabular}
\end{center}
\end{table}

%\newpage
%--------------------------------------------------------------------------------

\noindent
The following theorem provides evidence of
a new class of algebraic number fields with \(3\)-class group of type \(\mathrm{Cl}_3(K)\simeq (9,3)\)
whose \(3\)-class field tower consists of exactly three stages.

\begin{theorem}
\label{thm:ThreeStageTower93}
An imaginary quadratic field \(K=\mathbb{Q}(\sqrt{d})\)
with non-elementary \(3\)-class group \(\mathrm{Cl}_3(K)\simeq C_9\times C_3\) of rank two,
punctured capitulation type \(\mathrm{B}.18\), \(\varkappa(K)\sim (144;4)\),
and abelian type invariants \(\alpha_2(K)\) of second order of the shape in Formula
\eqref{eqn:SecondAQI21}
with either
\begin{equation}
\label{eqn:AQI1}
\tau_0=2^3, \ D_1=(32111)^3(221)^9, \ D_2=(321)^3(32)^9, \ D_3=(321)^3(32)^9, \ D_4=(321)^3(32)^9
\end{equation}
or
\begin{equation}
\label{eqn:AQI2}
\tau_0=2^3, \ D_1=(3211)^3(221)^9, \ D_2=(321)^3(32)^9, \ D_3=(321)^3(32)^9, \ D_4=(2211)^3(221)^9
\end{equation}
or
\begin{equation}
\label{eqn:AQI4}
\tau_0=2^3, \ D_1=(3211)^3(221)^9, \ D_2=(321)^3(32)^9, \ D_3=(321)^3(221)^9, \ D_4=(3111)^3(32)^9
\end{equation}
or
\begin{equation}
\label{eqn:AQI6}
\tau_0=321, \ D_1=(41111)^3(32)^9, \ D_2=(321)^3(32)^9, \ D_3=(321)^3(32)^9, \ D_4=(321)^3(32)^9
\end{equation}
possesses a finite \(3\)-class field tower
\[K=\mathrm{F}_3^0(K)<\mathrm{F}_3^1(K)<\mathrm{F}_3^2(K)<\mathrm{F}_3^3(K)=\mathrm{F}_3^\infty(K)\]
with precise length \(\ell_3(K)=3\).
\end{theorem}

%--------------------------------------------------------------------------------

%\noindent
In the following corollary,
Theorem
\ref{thm:ThreeStageTower93}
is supplemented by information on the Galois group
\(G=\mathrm{Gal}(\mathrm{F}_3^\infty(K)/K)\)
and its metabelianization
\(M=G/G^{\prime\prime}\simeq\mathrm{Gal}(\mathrm{F}_3^2(K)/K)\).

\begin{corollary}
\label{cor:ThreeStageTower93}
Let \(K\) be a field with properties as in the assumptions of Theorem
\ref{thm:ThreeStageTower93}.
Then the automorphism group \(G=\mathrm{Gal}(\mathrm{F}_3^\infty(K)/K)\)
of the full \(3\)-class field tower of \(K\) is a non-metabelian Schur \(\sigma\)-group
with soluble length \(\mathrm{sl}(G)=3\), order \(\#G=3^{21}\) and nilpotency class \(\mathrm{cl}(G)=9\).
The second \(3\)-class group \(M=\mathrm{Gal}(\mathrm{F}_3^2(K)/K)\) of \(K\)
is a metabelian \(\sigma\)-group of order \(\#M=3^{10}\) and nilpotency class \(\mathrm{cl}(M)=5\).
\end{corollary}

\begin{proof}
Formula
\eqref{eqn:AQI1}
leads to either \(\ell=2\), \(m=83\)
\cite[Lem. 10]{Ma2020b}
or \(\ell=7\), \(m=82\). \\
Formula
\eqref{eqn:AQI2}
leads to \(\ell\in\lbrace 4,6\rbrace\), \(m=83\) and \(108=81+27\) candidates for \(G\)
\cite[Lem. 6]{Ma2020b}. \\
Formula
\eqref{eqn:AQI4}
leads to either \(\ell\in\lbrace 10,13\rbrace\), \(m=82\)
or \(\ell\in\lbrace 12,15\rbrace\), \(m=83\)
\cite[Lem. 8]{Ma2020b}. \\
Formula
\eqref{eqn:AQI6}
leads to \(\ell=29\), \(m=85\) and \(27\) candidates for \(G\)
\cite[Lem. 10]{Ma2020b}. \\
Let \(B:=\langle 6561,165\rangle=\langle 729,10\rangle-\#2;2\) in the notation of
\cite{BEO2005,GNO2006}
be the common fork of the root paths of all finite \(3\)-groups \(G\)
with non-elementary bicyclic commutator quotient \(G/G^\prime\simeq C_9\times C_3\),
punctured transfer kernel type \(\mathrm{B}.18\),
\(\varkappa(G)\sim (144;4)\),
and logarithmic abelian quotient invariants of first order
\(\alpha_1(G)=\bigl(21;(\tau_0,211,211;211)\bigr)\)
with \(\tau_0\in\lbrace 222,321\rbrace\).
Then the candidates for \(G\) are given in the shape
\(B-\#4;\ell-\#2;k-\#4;j-\#1;i-\#2;h\)
with \(1\le\ell\le 72\), \(1\le k\le 41\), \(1\le j\le 27\), \(1\le i\le 5\),
where \(k\) is determined uniquely as a function \(k=k(\ell)\) of \(\ell\),
\(j\) runs through all possible values, 
\(i\) is determined uniquely as a function \(i=i(j)\) of \(j\),
and \(1\le h\le 3\)
\cite{Ma2020b}.
\end{proof}

%--------------------------------------------------------------------------------

\begin{example}
\label{exm:ThreeStageTower93}
The quadratic fields \(K\) with fundamental discriminants \\
\(d=-518\,835\) and \(\ell\in\lbrace 4,6\rbrace\), \\
respectively \(d\in\lbrace -553\,807, -763\,972, -876\,948\rbrace\) and \(\ell\in\lbrace 10,12,13,15\rbrace\), \\
respectively \(d=-150\,319\) and \(\ell=29\), \\
have punctured capitulation type \(\varkappa(K)\sim (144;4)\)
and are examples of field possessing a
\(3\)-class field tower with exactly three stages, \(\ell_3(K)=3\), of relative degrees
\[
\lbrack\mathrm{F}_3^3(K):\mathrm{F}_3^2(K)\rbrack=3^{11}, \quad
\lbrack\mathrm{F}_3^2(K):\mathrm{F}_3^1(K)\rbrack=3^7, \quad
\lbrack\mathrm{F}_3^1(K):\mathrm{F}_3^0(K)\rbrack=3^3,
\]
and Galois group \(\mathrm{Gal}(\mathrm{F}_3^\infty(K)/K)\) of order \(3^{21}\).
\end{example}

%--------------------------------------------------------------------------------

\begin{remark}
\label{rmk:ThreeStageTower93}
The quadratic fields \(K\) with fundamental discriminants \\
\(d=-761\,855\) and \(\ell\in\lbrace 3,5,11,14\rbrace\), \\
respectively \(d_K\in\lbrace -194\,703, -494\,771\rbrace\) and \(\ell\in\lbrace 16,19,23,25,27\rbrace\), \\
have punctured capitulation type \(\varkappa(K)\sim (144;4)\)
and \(3\le\ell_3(K)\le 4\).

The quadratic fields \(K\) with fundamental discriminants
\(d\in\lbrace -294\,983,-389\,435\rbrace\)
and punctured capitulation type \(\varkappa(K)\sim (144;4)\)
have an infinite \(3\)-class field tower.
\end{remark}

%\newpage
%--------------------------------------------------------------------------------

\section{\(3\)-groups with commutator quotient \((27,3)\)}
\label{s:31}

\noindent
In Table
\ref{tbl:31},
we list the second AQI \(\alpha_2\)
of the \(30\) non-metabelian step size-\(4\) descendants
\(F-\#4;\ell\) with \(43\le\ell\le 72\)
of the metabelian fork \(F=\langle 2187,3\rangle-\#2;10\).
The general structure of \(\alpha_2\) is
\begin{equation}
\label{eqn:2AQI31}
\alpha_2(G)=\lbrack 31;(421;32111,D_1),(311;32111,D_2),(311;32111,D_3);(221;32111,D_4)\rbrack,
\end{equation}
where each dodecuplet \(D_i\), \(1\le i\le 3\), consists of a triplet \(T_i^3\) and a nonet \(N_i^9\).
However, \(D_4\) consists of a triplet \(T_4^3\) and
either a nonet \(N_4^9\) or an octet \(O_4^8\) and a singlet \(S_4\).
The metabelianization \(M=G/G^{\prime\prime}\) is given by the step size-\(2\) descendant
\(F-\#2;m\) with \(m\in\lbrace 88,90\rbrace\) of the fork \(F\).
The smallest logarithmic order, soluble length, of a Schur \(\sigma\)-descendant \(S\) of \(G\)
is \(\mathrm{lo}(S)\), \(\mathrm{sl}(S)\).

%\newpage
%--------------------------------------------------------------------------------

\renewcommand{\arraystretch}{1.0}
\begin{table}[ht]
\caption{Invariants of \(G=\langle 2187,3\rangle-\#2;10-\#4;\ell\) with \(43\le\ell\le 72\)}
\label{tbl:31}
\begin{center}
\begin{tabular}{|r||l|l||l|l||l|l||l|l|l|l||c||c|c|}
\hline
 \(\ell\) & \(T_1\)   & \(N_1\) & \(T_2\)  & \(N_2\) & \(T_3\)  & \(N_3\) & \(T_4\)  & \(N_4\) & \(O_4\) & \(S_4\) & \(m\)  & \(\mathrm{lo}(S)\) & \(\mathrm{sl}(S)\) \\
\hline
   \(43\) & \(4211\)  & \(42\)  & \(421\)  & \(42\)  & \(421\)  & \(42\)  & \(3211\) & \(321\) &         &         & \(88\) & \(37\) & \(4\) \\
   \(58\) & \(4211\)  & \(42\)  & \(421\)  & \(42\)  & \(421\)  & \(42\)  & \(3211\) & \(321\) &         &         & \(90\) & \(37\) & \(4\) \\
\hline
   \(44\) & \(4211\)  & \(42\)  & \(421\)  & \(321\) & \(4111\) & \(42\)  & \(331\)  &         & \(321\) & \(222\) & \(88\) & \(22\) & \(3\) \\
   \(59\) & \(4211\)  & \(42\)  & \(421\)  & \(321\) & \(4111\) & \(42\)  & \(331\)  &         & \(321\) & \(222\) & \(90\) & \(22\) & \(3\) \\
\hline
   \(45\) & \(4211\)  & \(42\)  & \(421\)  & \(321\) & \(421\)  & \(321\) & \(3211\) & \(321\) &         &         & \(88\) & \(34\) & \(4\) \\
   \(62\) & \(4211\)  & \(42\)  & \(421\)  & \(321\) & \(421\)  & \(321\) & \(3211\) & \(321\) &         &         & \(90\) & \(34\) & \(4\) \\
\hline
   \(46\) & \(4211\)  & \(42\)  & \(421\)  & \(42\)  & \(4111\) & \(42\)  & \(331\)  &         & \(321\) & \(222\) & \(88\) & \(22\) & \(3\) \\
   \(50\) & \(4211\)  & \(42\)  & \(421\)  & \(42\)  & \(4111\) & \(42\)  & \(331\)  &         & \(321\) & \(222\) & \(88\) & \(22\) & \(3\) \\
   \(63\) & \(4211\)  & \(42\)  & \(421\)  & \(42\)  & \(4111\) & \(42\)  & \(331\)  &         & \(321\) & \(222\) & \(90\) & \(22\) & \(3\) \\
   \(65\) & \(4211\)  & \(42\)  & \(421\)  & \(42\)  & \(4111\) & \(42\)  & \(331\)  &         & \(321\) & \(222\) & \(90\) & \(22\) & \(3\) \\
\hline
   \(47\) & \(41111\) & \(42\)  & \(4111\) & \(42\)  & \(4111\) & \(42\)  & \(3211\) & \(321\) &         &         & \(88\) & \(\infty\) & \(\infty\) \\
   \(60\) & \(41111\) & \(42\)  & \(4111\) & \(42\)  & \(4111\) & \(42\)  & \(3211\) & \(321\) &         &         & \(90\) & \(\infty\) & \(\infty\) \\
\hline
   \(48\) & \(41111\) & \(42\)  & \(421\)  & \(42\)  & \(421\)  & \(321\) & \(331\)  &         & \(321\) & \(222\) & \(88\) & \(22\) & \(3\) \\
   \(61\) & \(41111\) & \(42\)  & \(421\)  & \(42\)  & \(421\)  & \(321\) & \(331\)  &         & \(321\) & \(222\) & \(90\) & \(22\) & \(3\) \\
\hline
   \(49\) & \(4211\)  & \(42\)  & \(421\)  & \(42\)  & \(3211\) & \(321\) & \(331\)  & \(321\) &         &         & \(88\) & \(25\) & \(3\) \\
   \(64\) & \(4211\)  & \(42\)  & \(421\)  & \(42\)  & \(3211\) & \(321\) & \(331\)  & \(321\) &         &         & \(90\) & \(25\) & \(3\) \\
\hline
   \(51\) & \(4211\)  & \(42\)  & \(421\)  & \(42\)  & \(421\)  & \(321\) & \(3211\) &         & \(321\) & \(222\) & \(88\) & \(22\) & \(3\) \\
   \(66\) & \(4211\)  & \(42\)  & \(421\)  & \(42\)  & \(421\)  & \(321\) & \(3211\) &         & \(321\) & \(222\) & \(90\) & \(22\) & \(3\) \\
\hline
   \(52\) & \(4211\)  & \(42\)  & \(421\)  & \(42\)  & \(4111\) & \(42\)  & \(331\)  & \(321\) &         &         & \(88\) & \(25\) & \(3\) \\
   \(70\) & \(4211\)  & \(42\)  & \(421\)  & \(42\)  & \(4111\) & \(42\)  & \(331\)  & \(321\) &         &         & \(90\) & \(25\) & \(3\) \\
\hline
   \(53\) & \(4211\)  & \(42\)  & \(421\)  & \(321\) & \(3211\) & \(321\) & \(331\)  &         & \(321\) & \(222\) & \(88\) & \(28\) & \(4\) \\
\hline
   \(54\) & \(4211\)  & \(42\)  & \(421\)  & \(42\)  & \(421\)  & \(42\)  & \(3211\) &         & \(321\) & \(222\) & \(88\) & \(22\) & \(3\) \\
   \(72\) & \(4211\)  & \(42\)  & \(421\)  & \(42\)  & \(421\)  & \(42\)  & \(3211\) &         & \(321\) & \(222\) & \(90\) & \(22\) & \(3\) \\
\hline
   \(55\) & \(41111\) & \(42\)  & \(421\)  & \(42\)  & \(421\)  & \(321\) & \(331\)  & \(321\) &         &         & \(88\) & \(25\) & \(3\) \\
   \(67\) & \(41111\) & \(42\)  & \(421\)  & \(42\)  & \(421\)  & \(321\) & \(331\)  & \(321\) &         &         & \(90\) & \(25\) & \(3\) \\
\hline
   \(56\) & \(41111\) & \(42\)  & \(421\)  & \(42\)  & \(421\)  & \(42\)  & \(331\)  &         & \(321\) & \(222\) & \(88\) & \(22\) & \(3\) \\
   \(68\) & \(41111\) & \(42\)  & \(421\)  & \(42\)  & \(421\)  & \(42\)  & \(331\)  &         & \(321\) & \(222\) & \(90\) & \(22\) & \(3\) \\
\hline
   \(57\) & \(41111\) & \(42\)  & \(4111\) & \(42\)  & \(3211\) & \(321\) & \(3211\) &         & \(321\) & \(222\) & \(88\) & \(\infty\) & \(\infty\) \\
   \(69\) & \(41111\) & \(42\)  & \(4111\) & \(42\)  & \(3211\) & \(321\) & \(3211\) &         & \(321\) & \(222\) & \(90\) & \(\infty\) & \(\infty\) \\
\hline
   \(71\) & \(4211\)  & \(42\)  & \(421\)  & \(321\) & \(3211\) & \(321\) & \(331\)  &         & \(321\) & \(222\) & \(90\) & \(25\) & \(3\) \\
\hline
\end{tabular}
\end{center}
\end{table}

%--------------------------------------------------------------------------------

\begin{theorem}
\label{thm:31}
The Schur \(\sigma\)-groups \(S\) with commutator quotient \(S/S^\prime\simeq (27,3)\),
punctured transfer kernel type \(\mathrm{B}.18\), \(\varkappa(S)\sim (144;4)\),
and first AQI \(\alpha_1(S)\sim (421,311,311;221)\)
are descendants of the \(30\) non-metabelian \(3\)-groups
\(G=\langle 2187,3\rangle-\#2;10-\#4;\ell\) whose invariants are listed in Table
\ref{tbl:31}.
In the case of finite order \(\mathrm{lo}(S)<\infty\),
their invariants coincide with those of the predecessor \(G\).
For \(22\le\mathrm{lo}(S)\le 25\) they have three stages \(\mathrm{sl}(S)=3\),
and for \(28\le\mathrm{lo}(S)<\infty\) they have four stages \(\mathrm{sl}(S)=4\).
For \(43\le\ell\le 57\) their metabelianization \(S/S^{\prime\prime}\simeq G/G^{\prime\prime}\)
is \(M=\langle 2187,3\rangle-\#2;10-\#2;88\),
and for \(58\le\ell\le 72\) it is \(M=\langle 2187,3\rangle-\#2;10-\#2;90\).
\end{theorem}

%--------------------------------------------------------------------------------

\noindent
The minimum \(\mathrm{lo}(S)=22\) occurs for \(14\) values \(\ell\),
\(25\) for \(7\), \(28\) for \(1\), \(34\) for \(2\), \(37\) for \(2\), and \(\infty\) for \(4\).

%\newpage
%--------------------------------------------------------------------------------

\section{Imaginary quadratic fields \(K\) with \(\mathrm{Cl}_3(K)\simeq C_{27}\times C_3\)}
\label{s:Imag31}

\noindent
The \(930\) imaginary quadratic fields \(K=\mathbb{Q}(\sqrt{d})\)
with fundamental discriminants \(-3\,000\,000<d<0\) and
\(3\)-class group \(\mathrm{Cl}_3(K)\simeq C_{27}\times C_3\)
were computed together with their punctured capitulation types \(\varkappa(K)\)
and first abelian type invariants \(\alpha_1(K)\)
by means of the computational algebra system Magma
\cite{MAGMA2021}
in \(19\,132\) seconds of CPU time, that is more than \(5\) hours.
In Table
\ref{tbl:Type31},
the first \(16\) cases with punctured capitulation type \(\mathrm{B}.18\),
\(\varkappa(K)\sim (144;4)\),
and \textit{uni-polarized} abelian type invariants \(\alpha_1(K)\) of first order
in the \textit{ground state} are listed.
Bi-polarized cases and excited states are excluded.

%\newpage
%--------------------------------------------------------------------------------

\renewcommand{\arraystretch}{1.1}
\begin{table}[ht]
\caption{Sixteen fields \(K=\mathbb{Q}(\sqrt{d})\) with \(\mathrm{Cl}_3(K)\simeq C_{27}\times C_3\) and \(\varkappa(K)\sim (144;4)\)}
\label{tbl:Type31}
\begin{center}
\begin{tabular}{|r|r|c||cc|}
\hline
 No.     & \(d\)            & factors          & \(\alpha_1(K)\)       & \\ %remark \\
\hline
  \(15\) &    \(-163\,736\) & \(2,97,211\)     & \((421,311,311;221)\) & \\
  \(25\) &    \(-218\,123\) & \(59,3\,697\)    & \((421,311,311;221)\) & \\
  \(75\) &    \(-428\,935\) & \(5,13,6\,599\)  & \((421,311,311;221)\) & \\
%\(111\) &    \(-566\,731\) & \(11,51\,521\)   & \((432,311,311;221)\) & first excited state \\
 \(121\) &    \(-615\,467\) & \(19,29,1\,117\) & \((421,311,311;221)\) & \\
 \(202\) &    \(-892\,459\) & \(31,28\,789\)   & \((421,311,311;221)\) & \\
 \(234\) &    \(-985\,727\) & \(463,2\,129\)   & \((421,311,311;221)\) & \\
 \(304\) & \(-1\,216\,407\) & \(3,47,8\,627\)  & \((421,311,311;221)\) & \\
 \(322\) & \(-1\,263\,279\) & \(3,421\,093\)   & \((421,311,311;221)\) & \\
 \(328\) & \(-1\,283\,531\) & \(701,1\,831\)   & \((421,311,311;221)\) & \\
 \(357\) & \(-1\,358\,087\) & prime            & \((421,311,311;221)\) & \\
 \(407\) & \(-1\,502\,187\) & \(3,500\,729\)   & \((421,311,311;221)\) & \\
 \(425\) & \(-1\,561\,043\) & \(11,191,743\)   & \((421,311,311;221)\) & \\
 \(433\) & \(-1\,588\,196\) & \(2,23,61,283\)  & \((421,311,311;221)\) & \\
 \(475\) & \(-1\,752\,787\) & \(67,26\,161\)   & \((421,311,311;221)\) & \\
 \(508\) & \(-1\,853\,828\) & \(2,463\,457\)   & \((421,311,311;221)\) & \\
 \(590\) & \(-2\,052\,195\) & \(3,5,136\,813\) & \((421,311,311;221)\) & \\
\hline
\end{tabular}
\end{center}
\end{table}

%\newpage
%--------------------------------------------------------------------------------

\noindent
In Table
\ref{tbl:SecondAQI31},
we give the abelian type invariants \(\alpha_2(K)\) of second order
of the \(16\) fields in the uni-polarized ground state
contained in Table
\ref{tbl:Type31}.
They were computed with the aid of Magma
\cite{MAGMA2021}
in \(74\,318\)  seconds of CPU time, that is nearly \(21\) hours.
The general structure of \(\alpha_2(K)\) is
\begin{equation}
\label{eqn:SecondAQI31}
\alpha_2(K)=\lbrack 31;(421;32111,D_1),(311;32111,D_2),(311;32111,D_3);(221;32111,D_4)\rbrack
\end{equation}
where each dodecuplet \(D_i\), \(1\le i\le 3\), consists of a triplet and a nonet,
and \(D_4\) consists of a triplet, an octet and a singlet.
A reference to Table
\ref{tbl:31}
is added. It usually admits the determination of
the length \(\ell_3(K)\) of the \(3\)-class field tower of \(K\).

%\newpage
%--------------------------------------------------------------------------------

\renewcommand{\arraystretch}{1.1}
\begin{table}[ht]
\caption{Details for the fields \(K=\mathbb{Q}(\sqrt{d})\) in Table \ref{tbl:Type31}}
\label{tbl:SecondAQI31}
\begin{center}
\begin{tabular}{|r||c|c|c|c||c|c|}
\hline
 No.     & \(D_1\)             & \(D_2\)            & \(D_3\)            & \(D_4\)                   & reference & \(\ell_3(K)\) \\
\hline
  \(15\) & \((42111)^3(42)^9\) & \((4211)^3(42)^9\) & \((3221)^3(321)^9\) & \((3311)^3(321)^8(222)\) & \(57,69\) var. & \(\infty\) \\
  \(25\) & \((4211)^3(42)^9\)  & \((421)^3(42)^9\)  & \((421)^3(321)^9\)  & \((3211)^3(321)^8(222)\) & \(51,66\) & \(3\) \\
  \(75\) & \((42111)^3(42)^9\) & \((4211)^3(42)^9\) & \((4211)^3(321)^9\) & \((3221)^3(321)^8(222)\) & \(57,69\) var. & \(\infty\) \\
 \(121\) & \((4211)^3(42)^9\)  & \((421)^3(321)^9\) & \((3211)^3(321)^9\) & \((331)^3(321)^8(222)\)  & \(53,71\) & \(4\) or \(3\) \\
 \(202\) & \((4211)^3(42)^9\)  & \((4111)^3(42)^9\) & \((421)^3(321)^9\)  & \((331)^3(321)^8(222)\)  & \(44,59\) & \(3\) \\
 \(234\) & \((42111)^3(42)^9\) & \((4211)^3(42)^9\) & \((4211)^3(321)^9\) & \((3311)^3(321)^8(222)\) & \(57,69\) var. & \(\infty\) \\
 \(304\) & \((51111)^3(42)^9\) & \((421)^3(42)^9\)  & \((421)^3(321)^9\)  & \((331)^3(321)^8(222)\)  & \(48,61\) var. & ? \\
 \(322\) & \((4211)^3(42)^9\)  & \((421)^3(42)^9\)  & \((4111)^3(42)^9\)  & \((331)^3(321)^9\)       & \(52,70\) & \(3\) \\
 \(328\) & \((51111)^3(42)^9\) & \((421)^3(42)^9\)  & \((421)^3(321)^9\)  & \((331)^3(321)^9\)       & \(55,67\) var. & ? \\
 \(357\) & \((4211)^3(42)^9\)  & \((421)^3(42)^9\)  & \((421)^3(321)^9\)  & \((3211)^3(321)^8(222)\) & \(51,66\) & \(3\) \\
 \(407\) & \((51111)^3(42)^9\) & \((421)^3(42)^9\)  & \((421)^3(321)^9\)  & \((331)^3(321)^8(222)\)  & \(48,61\) var. & ? \\
 \(425\) & \((42111)^3(42)^9\) & \((4211)^3(42)^9\) & \((4211)^3(321)^9\) & \((3311)^3(321)^8(222)\) & \(57,69\) var. & \(\infty\) \\
 \(433\) & \((4211)^3(42)^9\)  & \((421)^3(42)^9\)  & \((421)^3(42)^9\)   & \((3211)^3(321)^8(222)\) & \(43,58\) & \(4\) \\
 \(475\) & \((4211)^3(42)^9\)  & \((421)^3(42)^9\)  & \((4111)^3(42)^9\)  & \((331)^3(321)^8(222)\)  & \(46,50,63,65\) & \(3\) \\
 \(508\) & \((51111)^3(42)^9\) & \((421)^3(42)^9\)  & \((421)^3(321)^9\)  & \((331)^3(321)^9\)       & \(55,67\) var. & ? \\
 \(590\) & \((51111)^3(42)^9\) & \((421)^3(42)^9\)  & \((421)^3(42)^9\)   & \((331)^3(321)^8(222)\)  & \(56,68\) var. & ? \\
\hline
\end{tabular}
\end{center}
\end{table}

%--------------------------------------------------------------------------------

\begin{example}
\label{exm:Imaginary31a}
According to Tables
\ref{tbl:Type31}
and
\ref{tbl:SecondAQI31}
together with Theorem
\ref{thm:Imaginary31},
we get the following \(5\) examples of \(3\)-class field towers with precisely three stages,
\(\ell_3(K)=\mathrm{sl}(S)=3\):
\begin{itemize}
\item
\(d\in\lbrace -218\,123,-1\,358\,087\rbrace\) both with \(\ell\in\lbrace 51,66\rbrace\),
\item
\(d=-892\,459\) with \(\ell\in\lbrace 44,59\rbrace\),
\item
\(d=-1\,263\,279\) with \(\ell\in\lbrace 52,70\rbrace\) and \(\mathrm{lo}(S)=25\),
\item
\(d=-1\,752\,787\) with \(\ell\in\lbrace 46,50,63,65\rbrace\).
\end{itemize}
In contrast, the \(3\)-class field tower is infinite
for \(d\in\lbrace -163\,736,-428\,935,-985\,727,-1\,561\,043\rbrace\).
No statement is possible for
\(d\in\lbrace -1\,216\,407,-1\,283\,531,-1\,502\,187,-1\,853\,828,-2\,052\,195\rbrace\),
since the associated Schur \(\sigma\)-groups \(S\) are unknown.
\end{example}

%--------------------------------------------------------------------------------

\begin{example}
\label{exm:Imaginary31b}
As a particular highlight we point out the unique example
of a \(3\)-class field \textbf{tower with precisely four stages},
\(\ell_3(K)=\mathrm{sl}(S)=4\),
for \(d=-1\,588\,196\) with \(\ell\in\lbrace 43,58\rbrace\) and \(\mathrm{lo}(S)=37\).
As opposed, the precise length is unknown
for \(d=-615\,467\) with \(\ell\in\lbrace 53,71\rbrace\) and \(3\le\ell_3(K)=\mathrm{sl}(S)\le 4\).
\end{example}

%\newpage
%--------------------------------------------------------------------------------

\begin{theorem}
\label{thm:Imaginary31}
For an imaginary quadratic field \(K=\mathbb{Q}(\sqrt{d})\), \(d<0\),
with \(3\)-class group \(\mathrm{Cl}_3(K)\simeq (27,3)\)
and punctured capitulation type \(\mathrm{B}.18\), \(\varkappa\sim (144;4)\),
the \(3\)-class field tower consists of precisely three stages
with Schur \(\sigma\)-group \(G=\mathrm{Gal}(\mathrm{F}_3^\infty(K)/K)\)
of order \(\#G=3^{22}\) and nilpotency class \(\mathrm{cl}(G)=9\),
if the following conditions for the abelian type invariants \(\alpha_2(K)\) of second order in Formula
\eqref{eqn:SecondAQI31}
are satisfied.
In the notation of the SmallGroups database
\cite{BEO2005}
and the ANUPQ package
\cite{GNO2006},
the \(3\)-class field tower group is given by
\begin{equation}
\label{eqn:Schur31}
G\simeq\langle 2187,3\rangle-\#2;10-\#4;\ell-\#2;k(\ell)-\#4;j-\#1;i(j)-\#2;h,
\end{equation}
where \(43\le\ell\le 72\) is determined by the second AQI \(\alpha_2\),
\(1\le k\le 41\) is determined by \(\ell\),
\(1\le j\le 27\) is arbitrary,
\(1\le i\le 2\) is determined by \(j\),
and \(1\le h\le N\) is arbitrary below an upper bound \(N\in\lbrace 1,3\rbrace\) determined by \(\ell\).
\begin{itemize}
\item
\(\ell\in\lbrace 51,66\rbrace\), \(N=3\), i.e. \(162\) candidates for \(G\), \\
if \(D_1=(4211)^3(42)^9\), \(D_2=(421)^3(42)^9\), \(D_3=(421)^3(321)^9\), \(D_4=(3211)^3(321)^8(222)\);
\item
\(\ell\in\lbrace 44,59\rbrace\), \(N=1\), i.e. \(54\) candidates for \(G\), \\
if \(D_1=(4211)^3(42)^9\), \(D_2=(421)^3(321)^9\), \(D_3=(4111)^3(42)^9\), \(D_4=(331)^3(321)^8(222)\).
\end{itemize}
The metabelianization \(M=G/G^{\prime\prime}\simeq\mathrm{Gal}(\mathrm{F}_3^2(K)/K)\),
which is isomorphic to the second \(3\)-class group of \(K\),
has order \(\#M=3^{11}\), nilpotency class \(\mathrm{cl}(M)=5\) and is given by
\begin{equation}
\label{eqn:Sigma31}
M\simeq\langle 2187,3\rangle-\#2;10-\#2;m,
\end{equation}
where \(m=88\) if \(\ell\le 57\), and \(m=90\) if \(\ell\ge 58\).
\end{theorem}

\begin{proof}
Among the \(14\) descendants \(\langle 2187,3\rangle-\#2;10-\#4;\ell\)
which give rise to Schur \(\sigma\)-groups of minimal order \(3^{22}\),
that is \(\ell\in\lbrace 44,46,48,50,51,54,56,59,61,63,65,66,68,72\rbrace\),
the second AQI in the statements are unique.
It remains to check the other \(16\) values of \(43\le\ell\le 72\) with Tbl.
\ref{tbl:31}.
\end{proof}

%\newpage
%--------------------------------------------------------------------------------

\section{\(3\)-groups with commutator quotient \((81,3)\)}
\label{s:41}

\noindent
In Table
\ref{tbl:41},
we list the second AQI \(\alpha_2\)
of the \(30\) non-metabelian step size-\(4\) descendants
\(B-\#4;\ell\) with \(80\le\ell\le 109\)
of the metabelian fork \(B=\langle 2187,3\rangle-\#3;2\).
The general structure of \(\alpha_2\) is
\begin{equation}
\label{eqn:2AQI41}
\alpha_2(G)=\lbrack 41;(521;42111,D_1),(411;42111,D_2),(411;42111,D_3);(321;42111,D_4)\rbrack
\end{equation}
where each dodecuplet \(D_i\), \(1\le i\le 3\), consists of a triplet \(T_i^3\) and a nonet \(N_i^9\),
and \(D_4\) usually consists of a triplet \(T_4^3\), an octet \(O_4^8\) and a singlet \(S_4\).
The metabelianization \(M=G/G^{\prime\prime}\) is given by the step size-\(2\) descendant
\(B-\#2;m\) with \(m\in\lbrace 100,102\rbrace\) of the fork \(B\).
The smallest logarithmic order, soluble length, of a Schur \(\sigma\)-descendant \(S\) of \(G\)
is \(\mathrm{lo}(S)\), \(\mathrm{sl}(S)\).

%\newpage
%--------------------------------------------------------------------------------

\renewcommand{\arraystretch}{0.9}
\begin{table}[ht]
\caption{Invariants of \(G=\langle 2187,3\rangle-\#3;2-\#4;\ell\) with \(80\le\ell\le 109\)}
\label{tbl:41}
\begin{center}
\begin{tabular}{|r||l|l||l|l||l|l||l|l|l|l||c||c|c|}
\hline
 \(\ell\) & \(T_1\)   & \(N_1\) & \(T_2\)  & \(N_2\) & \(T_3\)  & \(N_3\) & \(T_4\)  & \(D_4\) & \(O_4\) & \(S_4\) & \(m\)   & \(\mathrm{lo}(S)\) & \(\mathrm{sl}(S)\) \\
\hline
   \(80\) & \(5211\)  & \(52\)  & \(521\)  & \(52\)  & \(521\)  & \(52\)  & \(3311\) & \((331)^6(322)^3\) & &      & \(100\) & \(46\) & \(4\) \\
   \(95\) & \(5211\)  & \(52\)  & \(521\)  & \(52\)  & \(521\)  & \(52\)  & \(3311\) & \((331)^6(322)^3\) & &      & \(102\) & \(46\) & \(4\) \\
\hline
   \(81\) & \(5211\)  & \(52\)  & \(521\)  & \(421\) & \(5111\) & \(52\)  & \(431\)  &         & \(421\) & \(322\) & \(100\) & \(23\) & \(3\) \\
   \(96\) & \(5211\)  & \(52\)  & \(521\)  & \(421\) & \(5111\) & \(52\)  & \(431\)  &         & \(421\) & \(322\) & \(102\) & \(23\) & \(3\) \\
\hline
   \(82\) & \(5211\)  & \(52\)  & \(521\)  & \(421\) & \(521\)  & \(421\) & \(3311\) & \((331)^6(322)^3\) & &      & \(100\) & \(40\) & \(4\) \\
   \(99\) & \(5211\)  & \(52\)  & \(521\)  & \(421\) & \(521\)  & \(421\) & \(3311\) & \((331)^6(322)^3\) & &      & \(102\) & \(40\) & \(4\) \\
\hline
   \(83\) & \(5211\)  & \(52\)  & \(521\)  & \(52\)  & \(5111\) & \(52\)  & \(431\)  &         & \(421\) & \(322\) & \(100\) & \(23\) & \(3\) \\
   \(87\) & \(5211\)  & \(52\)  & \(521\)  & \(52\)  & \(5111\) & \(52\)  & \(431\)  &         & \(421\) & \(322\) & \(100\) & \(23\) & \(3\) \\
  \(100\) & \(5211\)  & \(52\)  & \(521\)  & \(52\)  & \(5111\) & \(52\)  & \(431\)  &         & \(421\) & \(322\) & \(102\) & \(23\) & \(3\) \\
  \(102\) & \(5211\)  & \(52\)  & \(521\)  & \(52\)  & \(5111\) & \(52\)  & \(431\)  &         & \(421\) & \(322\) & \(102\) & \(23\) & \(3\) \\
\hline
   \(84\) & \(51111\) & \(52\)  & \(5111\) & \(52\)  & \(5111\) & \(52\)  & \(3311\) & \((331)^6(322)^3\) & &      & \(100\) & \(\infty\) & \(\infty\) \\
   \(97\) & \(51111\) & \(52\)  & \(5111\) & \(52\)  & \(5111\) & \(52\)  & \(3311\) & \((331)^6(322)^3\) & &      & \(102\) & \(\infty\) & \(\infty\) \\
\hline
   \(85\) & \(51111\) & \(52\)  & \(521\)  & \(52\)  & \(521\)  & \(421\) & \(431\)  &         & \(421\) & \(322\) & \(100\) & \(23\) & \(3\) \\
   \(98\) & \(51111\) & \(52\)  & \(521\)  & \(52\)  & \(521\)  & \(421\) & \(431\)  &         & \(421\) & \(322\) & \(102\) & \(23\) & \(3\) \\
\hline
   \(86\) & \(5211\)  & \(52\)  & \(521\)  & \(52\)  & \(4211\) & \(421\) & \(431\)  & \((331)^6(322)^3\) & &      & \(100\) & \(29\) & \(3\) \\
  \(101\) & \(5211\)  & \(52\)  & \(521\)  & \(52\)  & \(4211\) & \(421\) & \(431\)  & \((331)^6(322)^3\) & &      & \(102\) & \(29\) & \(3\) \\
\hline
   \(88\) & \(5211\)  & \(52\)  & \(521\)  & \(52\)  & \(521\)  & \(421\) & \(4211\) &         & \(421\) & \(322\) & \(100\) & \(23\) & \(3\) \\
  \(103\) & \(5211\)  & \(52\)  & \(521\)  & \(52\)  & \(521\)  & \(421\) & \(4211\) &         & \(421\) & \(322\) & \(102\) & \(23\) & \(3\) \\
\hline
   \(89\) & \(5211\)  & \(52\)  & \(521\)  & \(52\)  & \(5111\) & \(52\)  & \(431\)  & \((331)^6(322)^3\) & &      & \(100\) & \(29\) & \(3\) \\
  \(107\) & \(5211\)  & \(52\)  & \(521\)  & \(52\)  & \(5111\) & \(52\)  & \(431\)  & \((331)^6(322)^3\) & &      & \(102\) & \(29\) & \(3\) \\
\hline
   \(90\) & \(5211\)  & \(52\)  & \(521\)  & \(421\) & \(4211\) & \(421\) & \(431\)  &         & \(421\) & \(322\) & \(100\) & \(29\) & \(4\) \\
\hline
   \(91\) & \(5211\)  & \(52\)  & \(521\)  & \(52\)  & \(521\)  & \(52\)  & \(4211\) &         & \(421\) & \(322\) & \(100\) & \(23\) & \(3\) \\
  \(109\) & \(5211\)  & \(52\)  & \(521\)  & \(52\)  & \(521\)  & \(52\)  & \(4211\) &         & \(421\) & \(322\) & \(102\) & \(23\) & \(3\) \\
\hline
   \(92\) & \(51111\) & \(52\)  & \(521\)  & \(52\)  & \(521\)  & \(421\) & \(431\)  & \((331)^6(322)^3\) & &      & \(100\) & \(29\) & \(3\) \\
  \(104\) & \(51111\) & \(52\)  & \(521\)  & \(52\)  & \(521\)  & \(421\) & \(431\)  & \((331)^6(322)^3\) & &      & \(102\) & \(29\) & \(3\) \\
\hline
   \(93\) & \(51111\) & \(52\)  & \(521\)  & \(52\)  & \(521\)  & \(52\)  & \(431\)  &         & \(421\) & \(322\) & \(100\) & \(23\) & \(3\) \\
  \(105\) & \(51111\) & \(52\)  & \(521\)  & \(52\)  & \(521\)  & \(52\)  & \(431\)  &         & \(421\) & \(322\) & \(102\) & \(23\) & \(3\) \\
\hline
   \(94\) & \(51111\) & \(52\)  & \(5111\) & \(52\)  & \(4211\) & \(421\) & \(4211\) &         & \(421\) & \(322\) & \(100\) & \(\infty\) & \(\infty\) \\
  \(106\) & \(51111\) & \(52\)  & \(5111\) & \(52\)  & \(4211\) & \(421\) & \(4211\) &         & \(421\) & \(322\) & \(102\) & \(\infty\) & \(\infty\) \\
\hline
  \(108\) & \(5211\)  & \(52\)  & \(521\)  & \(421\) & \(4211\) & \(421\) & \(431\)  &         & \(421\) & \(322\) & \(102\) & \(26\) & \(3\) \\
\hline
\end{tabular}
\end{center}
\end{table}

%--------------------------------------------------------------------------------

\begin{theorem}
\label{thm:41}
The Schur \(\sigma\)-groups \(S\) with commutator quotient \(S/S^\prime\simeq (81,3)\),
punctured transfer kernel type \(\mathrm{B}.18\), \(\varkappa(S)\sim (144;4)\),
and first AQI \(\alpha_1(S)\sim (521,411,411;321)\)
are descendants of the \(30\) non-metabelian \(3\)-groups
\(G=\langle 2187,3\rangle-\#3;2-\#4;\ell\) whose invariants are listed in Table
\ref{tbl:41}.
In the case of finite order \(\mathrm{lo}(S)<\infty\),
their invariants coincide with those of the predecessor \(G\).
For \(23\le\mathrm{lo}(S)\le 26\) they have three stages \(\mathrm{sl}(S)=3\),
for \(40\le\mathrm{lo}(S)<\infty\) four stages \(\mathrm{sl}(S)=4\)
and for \(\mathrm{lo}(S)\le 29\) they have \(3\le\mathrm{sl}(S)\le 4\).
For \(80\le\ell\le 94\) their metabelianization \(S/S^{\prime\prime}\simeq G/G^{\prime\prime}\)
is \(M=\langle 2187,3\rangle-\#3;2-\#2;100\),
and for \(95\le\ell\le 109\) it is \(M=\langle 2187,3\rangle-\#3;2-\#2;102\).
\end{theorem}

\noindent
The minimum \(\mathrm{lo}(S)=23\) occurs for \(14\) values \(\ell\),
\(26\) for \(1\), \(29\) for \(7\), \(40\) for \(2\), \(46\) for \(2\), and \(\infty\) for \(4\).

%--------------------------------------------------------------------------------

\begin{remark}
\label{rmk:31And41}
Table
\ref{tbl:31}
and Theorem
\ref{thm:31}
were completed on \(24\) August \(2021\).
After the discovery of the fork \(B=\langle 2187,3\rangle-\#3;2\)
on \(26\) August \(2021\), Table
\ref{tbl:41}
could be computed immediately:
For all \(1\le e\le 4\),
exemplary representatives of multiplets of Schur \(\sigma\)-groups \(S\)
with commutator quotient \(S/S^\prime\simeq (3^e,3)\)
can be found according to the \textit{principle of extremal root paths}
(see Figure
\ref{fig:SchurSigma}),
\[
S\buildrel s=2 \over\longrightarrow
\pi(S)\buildrel s=1 \over\longrightarrow
\pi^2(S)\buildrel s=4 \over\longrightarrow
\pi^3(S)\buildrel s=2 \over\longrightarrow
\pi^4(S)\buildrel s=4 \over\longrightarrow
\pi^5(S)=B.
\]
This is not possible any longer for \(e\ge 5\),
due to the beginning discrepancy between parents \(\pi(G)\) and \(p\)-parents \(\pi_p(G)\)
of finite \(3\)-groups \(G\).
\end{remark}

%\newpage
%--------------------------------------------------------------------------------

\section{Imaginary quadratic fields \(K\) with \(\mathrm{Cl}_3(K)\simeq C_{81}\times C_3\)}
\label{s:Imag41}

\noindent
The \(2174\) imaginary quadratic fields \(K=\mathbb{Q}(\sqrt{d})\)
with fundamental discriminants \(-20\,000\,000<d<0\) and
\(3\)-class group \(\mathrm{Cl}_3(K)\simeq C_{81}\times C_3\)
were computed by means of the computational algebra system Magma
\cite{MAGMA2021}
in \(141\,586\) seconds of CPU time, that is nearly two days.
In Table
\ref{tbl:Type41},
the first eight cases with punctured capitulation type \(\mathrm{B}.18\),
\(\varkappa(K)\sim (144;4)\),
are listed.
The abelian quotient invariants \(\alpha_1(K)\) of first order
of only six of them are \textit{uni-polarized} and in the \textit{ground state}.

%\newpage
%--------------------------------------------------------------------------------

\renewcommand{\arraystretch}{1.1}
\begin{table}[ht]
\caption{Eight fields \(K=\mathbb{Q}(\sqrt{d})\) with \(\mathrm{Cl}_3(K)\simeq C_{81}\times C_3\) and \(\varkappa(K)\sim (144;4)\)}
\label{tbl:Type41}
\begin{center}
\begin{tabular}{|r|r|c||c|c|}
\hline
 No.     & \(d\)             & factors        & \(\alpha_1(K)\)       & remark \\
\hline
  \(31\) &    \(-936\,311\) & prime           & \((521,411,411;321)\) & \\
  \(49\) & \(-1\,240\,879\) & \(107,11\,597\) & \((521,411,411;321)\) & \\
  \(57\) & \(-1\,437\,179\) & \(19,75\,641\)  & \((521,411,411;321)\) & \\
  \(78\) & \(-1\,723\,864\) & \(2,215\,483\)  & \((521,411,411;321)\) & \\
  \(80\) & \(-1\,749\,655\) & \(5,349\,931\)  & \((521,411,411;321)\) & \\
  \(86\) & \(-1\,818\,223\) & \(11,165\,293\) & \((521,411,411;321)\) & \\
  \(96\) & \(-1\,854\,319\) & \(281,6\,599\)  & \((532,411,411;321)\) & first excited state \\
 \(109\) & \(-2\,003\,179\) & \(61,32\,839\)  & \((521,411,411;332)\) & bi-polarized \\
\hline
\end{tabular}
\end{center}
\end{table}

%\newpage
%--------------------------------------------------------------------------------

\noindent
In Table
\ref{tbl:SecondAQI41},
we give the abelian quotient invariants \(\alpha_2(K)\) of second order
of the six fields in the uni-polarized ground state
contained in Table
\ref{tbl:Type41}.
The general structure of \(\alpha_2(K)\) is the following
\begin{equation}
\label{eqn:SecondAQI41}
\alpha_2(K)=\lbrack 41;(521;42111,D_1),(411;42111,D_2),(411;42111,D_3);(321;42111,D_4)\rbrack
\end{equation}
where each dodecuplet \(D_i\), \(1\le i\le 3\), consists of a triplet and a nonet,
and \(D_4\) usually consists of a triplet, an octet and a singlet.
But the constitution of \(D_4\) may occasionally be irregular.

%\newpage
%--------------------------------------------------------------------------------

\renewcommand{\arraystretch}{1.1}
\begin{table}[ht]
\caption{Details for six fields \(K=\mathbb{Q}(\sqrt{d})\) in Table \ref{tbl:Type41}}
\label{tbl:SecondAQI41}
\begin{center}
\begin{tabular}{|r||c|c|c|c|c|}
\hline
 No.    & \(D_1\)             & \(D_2\)            & \(D_3\)            & \(D_4\)                   & remark \\
\hline
 \(31\) & \((5211)^3(52)^9\)  & \((5111)^3(52)^9\) & \((521)^3(421)^9\) & \((431)^3(421)^8(322)\)   & ref. \(81,96\) \\
 \(49\) & \((5211)^3(52)^9\)  & \((5111)^3(52)^9\) & \((521)^3(52)^9\)  & \((431)^3(421)^8(322)\)   & ref. \(83,87,100,102\) \\
 \(57\) & \((5211)^3(52)^9\)  & \((521)^3(52)^9\)  & \((521)^3(421)^9\) & \((4211)^3(421)^8(322)\)  & ref. \(88,103\) \\
 \(78\) & \(\)                & \(\)               & \(\)               & \((431)^3(421)^8(322)\)   & Magma int. err. \\
 \(80\) & \((61111)^3(52)^9\) & \((521)^3(52)^9\)  & \((521)^3(421)^9\) & \((431)^3(331)^6(322)^3\) & irregular \\
 \(86\) & \((5211)^3(52)^9\)  & \((5111)^3(52)^9\) & \((521)^3(52)^9\)  & \((431)^3(421)^8(322)\)   & ref. \(83,87,100,102\) \\
\hline
\end{tabular}
\end{center}
\end{table}

%\newpage
%--------------------------------------------------------------------------------

\begin{theorem}
\label{thm:Imaginary41}
For an imaginary quadratic field \(K=\mathbb{Q}(\sqrt{d})\), \(d<0\),
with \(3\)-class group \(\mathrm{Cl}_3(K)\simeq (81,3)\)
and punctured capitulation type \(\mathrm{B}.18\), \(\varkappa\sim (144;4)\),
the \(3\)-class field tower consists of precisely three stages
with Schur \(\sigma\)-group \(G=\mathrm{Gal}(\mathrm{F}_3^\infty(K)/K)\)
of order \(\#G=3^{23}\) and nilpotency class \(\mathrm{cl}(G)=9\),
if the following conditions for the abelian quotient invariants \(\alpha_2(K)\) of second order in Formula
\eqref{eqn:SecondAQI41}
are satisfied.
In the notation of the SmallGroups database
\cite{BEO2005}
and the ANUPQ package
\cite{GNO2006},
the \(3\)-class field tower group is given by
\begin{equation}
\label{eqn:Schur41}
G\simeq\langle 2187,3\rangle-\#3;2-\#4;\ell-\#2;k(\ell)-\#4;j-\#1;i(j)-\#2;h,
\end{equation}
where \(80\le\ell\le 109\) is determined by the second AQI \(\alpha_2\),
\(1\le k\le 41\) is determined by \(\ell\),
\(1\le j\le 27\) is arbitrary,
\(1\le i\le 2\) is determined by \(j\),
and \(1\le h\le N\) is arbitrary below an upper bound \(N\in\lbrace 1,3\rbrace\) determined by \(\ell\).
\begin{itemize}
\item
\(\ell\in\lbrace 81,96\rbrace\), \(N=1\), i.e. \(54\) candidates for \(G\), \\
if \(D_1=(5211)^3(52)^9\), \(D_2=(521)^3(421)^9\), \(D_3=(5111)^3(52)^9\), \(D_4=(431)^3(421)^8(322)\);
\item
\(\ell\in\lbrace 83,87,100,102\rbrace\), \(N=3\), i.e. \(324\) candidates for \(G\), \\
if \(D_1=(5211)^3(52)^9\), \(D_2=(521)^3(52)^9\), \(D_3=(5111)^3(52)^9\), \(D_4=(431)^3(421)^8(322)\);
\item
\(\ell\in\lbrace 88,103\rbrace\), \(N=3\), i.e. \(162\) candidates for \(G\), \\
if \(D_1=(5211)^3(52)^9\), \(D_2=(521)^3(52)^9\), \(D_3=(521)^3(421)^9\), \(D_4=(4211)^3(421)^8(322)\).
\end{itemize}
The metabelianization \(M=G/G^{\prime\prime}\simeq\mathrm{Gal}(\mathrm{F}_3^2(K)/K)\),
which is isomorphic to the second \(3\)-class group of \(K\),
has order \(\#M=3^{12}\), nilpotency class \(\mathrm{cl}(M)=5\) and is given by
\begin{equation}
\label{eqn:Sigma41}
M\simeq\langle 2187,3\rangle-\#3;2-\#2;m,
\end{equation}
where \(m=100\) if \(\ell\le 93\), and \(m=102\) if \(\ell\ge 96\).
\end{theorem}

\begin{proof}
The root \(\langle 2187,3\rangle-\#3;2\) can be viewed as
usual descendant of \(\langle 6561,216\rangle\) with step size \(s=2\).
Among the \(14\) descendants \(\langle 2187,3\rangle-\#3;2-\#4;\ell\)
which give rise to Schur \(\sigma\)-groups of minimal order \(3^{23}\),
that is \(\ell\in\lbrace 81,83,85,87,88,91,93,96,98,100,102,103,105,109\rbrace\),
the second AQI in the statements are unique.
It remains to check the other \(16\) values of \(80\le\ell\le 109\).
\end{proof}

%--------------------------------------------------------------------------------

\begin{example}
\label{exm:Imaginary41}
According to Tables
\ref{tbl:Type41}
and
\ref{tbl:SecondAQI41}
together with Theorem
\ref{thm:Imaginary41},
we get the following \(4\) examples of \(3\)-class field towers with precisely three stages,
\(\ell_3(K)=\mathrm{sl}(S)=3\):
\begin{itemize}
\item
\(d=-936\,311\) with \(\ell\in\lbrace 81,96\rbrace\),
\item
\(d=-1\,240\,879\) and \(d=-1\,818\,223\) both with \(\ell\in\lbrace 83,87,100,102\rbrace\),
\item
\(d=-1\,437\,179\) with \(\ell\in\lbrace 88,103\rbrace\).
\end{itemize}
In contrast, no statement is possible for \(d=-1\,749\,655\).
\end{example}

%\newpage
%--------------------------------------------------------------------------------

\section{Motivation for seeking the new periodicities of Schur \(\sigma\)-groups}
\label{s:Periodicities}

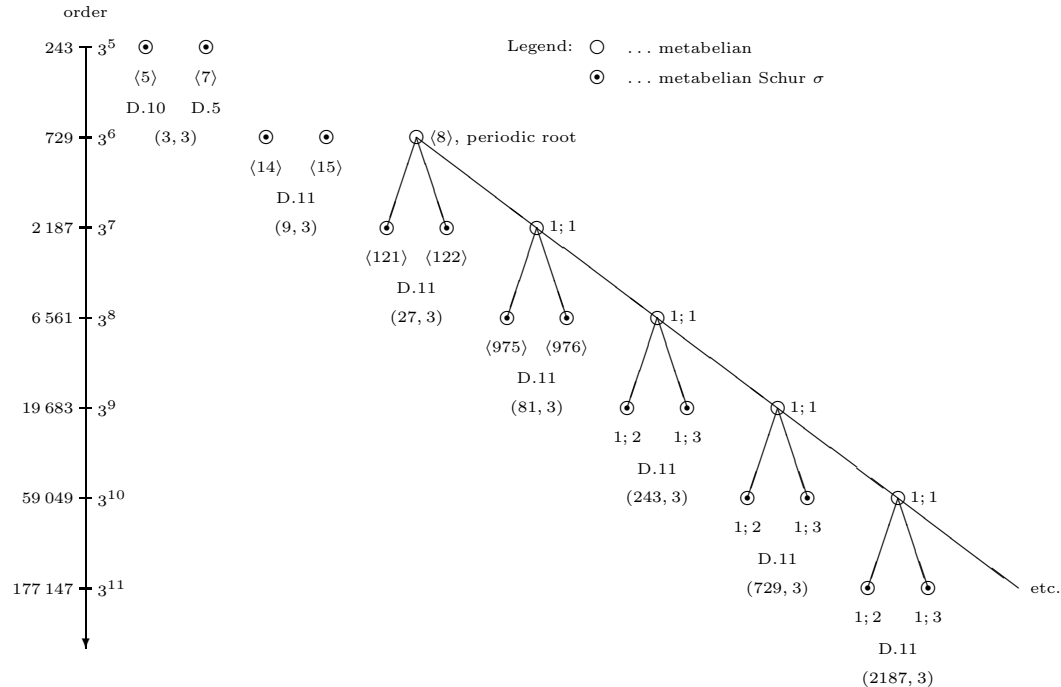
\begin{figure}[hb]
\caption{Periodic metabelian Schur \(\sigma\)-groups \(G\) with \(G/G^\prime\simeq (3^e,3)\), \(e\ge 3\)}
\label{fig:SchurSigma3}

{\tiny

\setlength{\unitlength}{0.8cm}
\begin{picture}(16,11)(-10,-8)

% scale of orders
\put(-10,2.5){\makebox(0,0)[cb]{order}}

\put(-10,2){\line(0,-1){9}}
\multiput(-10.1,2)(0,-1.5){7}{\line(1,0){0.2}}

%***********************************************
\put(-10.2,2){\makebox(0,0)[rc]{\(243\)}}
\put(-9.8,2){\makebox(0,0)[lc]{\(3^5\)}}
\put(-10.2,0.5){\makebox(0,0)[rc]{\(729\)}}
\put(-9.8,0.5){\makebox(0,0)[lc]{\(3^6\)}}
\put(-10.2,-1){\makebox(0,0)[rc]{\(2\,187\)}}
\put(-9.8,-1){\makebox(0,0)[lc]{\(3^7\)}}
\put(-10.2,-2.5){\makebox(0,0)[rc]{\(6\,561\)}}
\put(-9.8,-2.5){\makebox(0,0)[lc]{\(3^8\)}}
\put(-10.2,-4){\makebox(0,0)[rc]{\(19\,683\)}}
\put(-9.8,-4){\makebox(0,0)[lc]{\(3^9\)}}
\put(-10.2,-5.5){\makebox(0,0)[rc]{\(59\,049\)}}
\put(-9.8,-5.5){\makebox(0,0)[lc]{\(3^{10}\)}}
\put(-10.2,-7){\makebox(0,0)[rc]{\(177\,147\)}}
\put(-9.8,-7){\makebox(0,0)[lc]{\(3^{11}\)}}
%***********************************************

\put(-10,-7){\vector(0,-1){1}}

% metabelian vertices
\put(-9,2){\circle{0.2}}
\put(-9,2){\circle*{0.1}}
\put(-8,2){\circle{0.2}}
\put(-8,2){\circle*{0.1}}

\put(-7,0.5){\circle{0.2}}
\put(-7,0.5){\circle*{0.1}}
\put(-6,0.5){\circle{0.2}}
\put(-6,0.5){\circle*{0.1}}

\put(-4.5,0.5){\circle{0.2}}
\put(-5,-1){\circle{0.2}}
\put(-5,-1){\circle*{0.1}}
\put(-4,-1){\circle{0.2}}
\put(-4,-1){\circle*{0.1}}

\put(-2.5,-1){\circle{0.2}}
\put(-3,-2.5){\circle{0.2}}
\put(-3,-2.5){\circle*{0.1}}
\put(-2,-2.5){\circle{0.2}}
\put(-2,-2.5){\circle*{0.1}}

\put(-0.5,-2.5){\circle{0.2}}
\put(-1,-4){\circle{0.2}}
\put(-1,-4){\circle*{0.1}}
\put(0,-4){\circle{0.2}}
\put(0,-4){\circle*{0.1}}

\put(1.5,-4){\circle{0.2}}
\put(1,-5.5){\circle{0.2}}
\put(1,-5.5){\circle*{0.1}}
\put(2,-5.5){\circle{0.2}}
\put(2,-5.5){\circle*{0.1}}

\put(3.5,-5.5){\circle{0.2}}
\put(3,-7){\circle{0.2}}
\put(3,-7){\circle*{0.1}}
\put(4,-7){\circle{0.2}}
\put(4,-7){\circle*{0.1}}

%***********************************************
% directed edges 

\put(-4.5,0.5){\line(-1,-3){0.5}}
\put(-4.5,0.5){\line(1,-3){0.5}}
\put(-4.5,0.5){\line(4,-3){2}}

\put(-2.5,-1){\line(-1,-3){0.5}}
\put(-2.5,-1){\line(1,-3){0.5}}
\put(-2.5,-1){\line(4,-3){2}}

\put(-0.5,-2.5){\line(-1,-3){0.5}}
\put(-0.5,-2.5){\line(1,-3){0.5}}
\put(-0.5,-2.5){\line(4,-3){2}}

\put(1.5,-4){\line(-1,-3){0.5}}
\put(1.5,-4){\line(1,-3){0.5}}
\put(1.5,-4){\line(4,-3){2}}

\put(3.5,-5.5){\line(-1,-3){0.5}}
\put(3.5,-5.5){\line(1,-3){0.5}}
\put(3.5,-5.5){\line(4,-3){2}}
\put(5.7,-7){\makebox(0,0)[lc]{etc.}}

%***********************************************
% identifiers in the SmallGroups Database

% metabelian vertices
\put(-9,1.5){\makebox(0,0)[cc]{\(\langle 5\rangle\)}}
\put(-8,1.5){\makebox(0,0)[cc]{\(\langle 7\rangle\)}}
\put(-9,1){\makebox(0,0)[cc]{\(\mathrm{D}.10\)}}
\put(-8,1){\makebox(0,0)[cc]{\(\mathrm{D}.5\)}}
\put(-8.5,0.5){\makebox(0,0)[cc]{\((3,3)\)}}

\put(-7,0){\makebox(0,0)[cc]{\(\langle 14\rangle\)}}
\put(-6,0){\makebox(0,0)[cc]{\(\langle 15\rangle\)}}
\put(-6.5,-0.5){\makebox(0,0)[cc]{\(\mathrm{D}.11\)}}
\put(-6.5,-1){\makebox(0,0)[cc]{\((9,3)\)}}

\put(-4.3,0.5){\makebox(0,0)[lc]{\(\langle 8\rangle\), periodic root}}
\put(-5,-1.5){\makebox(0,0)[cc]{\(\langle 121\rangle\)}}
\put(-4,-1.5){\makebox(0,0)[cc]{\(\langle 122\rangle\)}}
\put(-4.5,-2){\makebox(0,0)[cc]{\(\mathrm{D}.11\)}}
\put(-4.5,-2.5){\makebox(0,0)[cc]{\((27,3)\)}}

\put(-2.3,-1){\makebox(0,0)[lc]{\(1;1\)}}
\put(-3,-3){\makebox(0,0)[cc]{\(\langle 975\rangle\)}}
\put(-2,-3){\makebox(0,0)[cc]{\(\langle 976\rangle\)}}
\put(-2.5,-3.5){\makebox(0,0)[cc]{\(\mathrm{D}.11\)}}
\put(-2.5,-4){\makebox(0,0)[cc]{\((81,3)\)}}

\put(-0.3,-2.5){\makebox(0,0)[lc]{\(1;1\)}}
\put(-1,-4.5){\makebox(0,0)[cc]{\(1;2\)}}
\put(0,-4.5){\makebox(0,0)[cc]{\(1;3\)}}
\put(-0.5,-5){\makebox(0,0)[cc]{\(\mathrm{D}.11\)}}
\put(-0.5,-5.5){\makebox(0,0)[cc]{\((243,3)\)}}

\put(1.7,-4){\makebox(0,0)[lc]{\(1;1\)}}
\put(1,-6){\makebox(0,0)[cc]{\(1;2\)}}
\put(2,-6){\makebox(0,0)[cc]{\(1;3\)}}
\put(1.5,-6.5){\makebox(0,0)[cc]{\(\mathrm{D}.11\)}}
\put(1.5,-7){\makebox(0,0)[cc]{\((729,3)\)}}

\put(3.7,-5.5){\makebox(0,0)[lc]{\(1;1\)}}
\put(3,-7.5){\makebox(0,0)[cc]{\(1;2\)}}
\put(4,-7.5){\makebox(0,0)[cc]{\(1;3\)}}
\put(3.5,-8){\makebox(0,0)[cc]{\(\mathrm{D}.11\)}}
\put(3.5,-8.5){\makebox(0,0)[cc]{\((2187,3)\)}}

%***********************************************
% legend

\put(-3,2){\makebox(0,0)[lc]{Legend:}}

\put(-1.5,2){\circle{0.2}}
\put(-1,2){\makebox(0,0)[lc]{\(\ldots\) metabelian}}

\put(-1.5,1.5){\circle{0.2}}
\put(-1.5,1.5){\circle*{0.1}}
\put(-1,1.5){\makebox(0,0)[lc]{\(\ldots\) metabelian Schur \(\sigma\)}}

\end{picture}

}

\end{figure}

\noindent
In our previous work
\cite[\S\ 7, Thm. 4 and Thm. 7]{Ma2021},
we found a periodicity of pairs of \textit{metabelian} Schur \(\sigma\)-groups \(G\)
with \(G/G^\prime\simeq (3^e,3)\), \(e\ge 3\),
and type \(\mathrm{D}.11\), \(\varkappa\sim (124;1)\),
which is illustrated by Figure
\ref{fig:SchurSigma3}.

%\newpage
%--------------------------------------------------------------------------------

%\newpage
%--------------------------------------------------------------------------------

In the Figures
\ref{fig:SchurSigma3} --
%\ref{fig:SchurSigma5},
%\ref{fig:SchurSigma}, and
\ref{fig:SchurSigma59},
all directed edges lead from descendants \(D\) to \(p\)-parents \(\pi_p(D)=D/P_{c_p-1}(D)\),
rather than to parents \(\pi(D)=D/\gamma_c(D)\).
The figures admit actual descendant construction.

\begin{figure}[ht]
\caption{Schur \(\sigma\)-groups \(G\) with \(\varrho(G)\sim (2,2,3;3)\), \(G/G^\prime\simeq (3^e,3)\), \(2\le e\le 7\)}
\label{fig:SchurSigma5}

{\tiny

\setlength{\unitlength}{0.9cm}
\begin{picture}(16,19.6)(-11,-16.5)

% scale of orders
\put(-10,2.5){\makebox(0,0)[cb]{order}}

\put(-10,2){\line(0,-1){17}}
\multiput(-10.1,2)(0,-1.5){13}{\line(1,0){0.2}}

%***********************************************
\put(-10.2,2){\makebox(0,0)[rc]{\(9\)}}
\put(-9.8,2){\makebox(0,0)[lc]{\(3^2\)}}
\put(-10.2,0.5){\makebox(0,0)[rc]{\(27\)}}
\put(-9.8,0.5){\makebox(0,0)[lc]{\(3^3\)}}
\put(-10.2,-1){\makebox(0,0)[rc]{\(81\)}}
\put(-9.8,-1){\makebox(0,0)[lc]{\(3^4\)}}
\put(-10.2,-2.5){\makebox(0,0)[rc]{\(243\)}}
\put(-9.8,-2.5){\makebox(0,0)[lc]{\(3^5\)}}
\put(-10.2,-4){\makebox(0,0)[rc]{\(729\)}}
\put(-9.8,-4){\makebox(0,0)[lc]{\(3^6\)}}
\put(-10.2,-5.5){\makebox(0,0)[rc]{\(2\,187\)}}
\put(-9.8,-5.5){\makebox(0,0)[lc]{\(3^7\)}}
\put(-10.2,-7){\makebox(0,0)[rc]{\(6\,561\)}}
\put(-9.8,-7){\makebox(0,0)[lc]{\(3^8\)}}
\put(-10.2,-8.5){\makebox(0,0)[rc]{\(19\,683\)}}
\put(-9.8,-8.5){\makebox(0,0)[lc]{\(3^9\)}}
\put(-10.2,-10){\makebox(0,0)[rc]{\(59\,049\)}}
\put(-9.8,-10){\makebox(0,0)[lc]{\(3^{10}\)}}
\put(-10.2,-11.5){\makebox(0,0)[rc]{\(177\,147\)}}
\put(-9.8,-11.5){\makebox(0,0)[lc]{\(3^{11}\)}}
\put(-10.2,-13){\makebox(0,0)[rc]{\(531\,441\)}}
\put(-9.8,-13){\makebox(0,0)[lc]{\(3^{12}\)}}
\put(-10.2,-14.5){\makebox(0,0)[rc]{\(1\,594\,323\)}}
\put(-9.8,-14.5){\makebox(0,0)[lc]{\(3^{13}\)}}
\put(-10.2,-16){\makebox(0,0)[rc]{\(4\,782\,969\)}}
\put(-9.8,-16){\makebox(0,0)[lc]{\(3^{14}\)}}
%***********************************************
% vertices

\put(-10,-15){\vector(0,-1){2}}

% abelian root vertex
\put(-9.1,1.9){\framebox(0.2,0.2){}}
\put(-9,2){\circle*{0.1}}

% metabelian vertices
\put(-9,0.5){\circle{0.2}}
\put(-9,-2.5){\circle{0.2}}
\put(-9,-4){\circle*{0.2}}
\put(-8.5,-5.5){\circle{0.2}}

\put(-7,-1){\circle{0.2}}
\put(-7,-4){\circle{0.2}}
\put(-7,-5.5){\circle*{0.2}}
\put(-6.5,-7){\circle{0.2}}

\put(-5,-4){\circle{0.2}}
\put(-5,-7){\circle*{0.2}}
\put(-4.5,-8.5){\circle{0.2}}

\put(-3,-7){\circle{0.2}}
\put(-3,-10){\circle{0.2}}

\put(-1,-10){\circle{0.2}}
\put(-1,-11.5){\circle{0.2}}

\put(1,-11.5){\circle{0.2}}
\put(1,-13){\circle{0.2}}

\put(3,-13){\circle{0.2}}
\put(3,-14.5){\circle{0.2}}

% non-metabelian vertices
\put(-9.1,-7.1){\framebox(0.2,0.2){}}
\put(-7.1,-8.6){\framebox(0.2,0.2){}}
\put(-5.1,-10.1){\framebox(0.2,0.2){}}
\put(-3.1,-11.6){\framebox(0.2,0.2){}}
\put(-1.1,-13.1){\framebox(0.2,0.2){}}
\put(0.9,-14.6){\framebox(0.2,0.2){}}
\put(2.9,-16.1){\framebox(0.2,0.2){}}

%***********************************************
% directed edges 

\put(-9,2){\line(0,-1){1.5}}
\put(-9,0.5){\line(0,-1){3}}
\put(-9,-2.5){\line(0,-1){1.5}}
\put(-9,-4){\line(0,-1){3}}
\put(-9,-4){\line(1,-3){0.5}}

\put(-9,2){\line(2,-3){2}}
\put(-7,-1){\line(0,-1){3}}
\put(-7,-4){\line(0,-1){1.5}}
\put(-7,-5.5){\line(0,-1){3}}
\put(-7,-5.5){\line(1,-3){0.5}}

\put(-7,-1){\line(2,-3){2}}
\put(-5,-4){\line(0,-1){3}}
\put(-5,-7){\line(0,-1){3}}
\put(-5,-7){\line(1,-3){0.5}}

\put(-5,-4){\line(2,-3){2}}
\put(-3,-7){\line(0,-1){3}}
\put(-3,-10){\line(0,-1){1.5}}

\put(-3,-7){\line(2,-3){2}}
\put(-1,-10){\line(0,-1){1.5}}
\put(-1,-11.5){\line(0,-1){1.5}}

\put(-1,-10){\line(4,-3){2}}
\put(1,-11.5){\line(0,-1){1.5}}
\put(1,-13){\line(0,-1){1.5}}

\put(1,-11.5){\line(4,-3){2}}
\put(3,-13){\line(0,-1){1.5}}
\put(3,-14.5){\line(0,-1){1.5}}

\put(3,-13){\line(4,-3){1}}

%***********************************************
% identifiers in the SmallGroups Database

% abelian root vertex
\put(-8.7,2){\makebox(0,0)[lc]{\(\langle 2\rangle\)}}

% metabelian vertices
\put(-8.7,0.5){\makebox(0,0)[lc]{\(\langle 3\rangle\)}}
\put(-8.7,-2.5){\makebox(0,0)[lc]{\(\langle 8\rangle\)}}
\put(-8.7,-4){\makebox(0,0)[lc]{\(\langle 54\rangle\)}}
\put(-8.2,-5.5){\makebox(0,0)[lc]{\(\langle 304\rangle\)}}

\put(-6.7,-1){\makebox(0,0)[lc]{\(\langle 3\rangle\)}}
\put(-6.7,-4){\makebox(0,0)[lc]{\(\langle 13\rangle\)}}
\put(-6.7,-5.5){\makebox(0,0)[lc]{\(\langle 168\rangle\)}}
\put(-6.2,-7){\makebox(0,0)[lc]{\(\langle 1689\rangle\)}}

\put(-4.7,-4){\makebox(0,0)[lc]{\(\langle 7\rangle\)}}
\put(-4.7,-7){\makebox(0,0)[lc]{\(\langle 98\rangle\)}}
\put(-4.2,-8.5){\makebox(0,0)[lc]{\(1;4\)}}

\put(-2.7,-7){\makebox(0,0)[lc]{\(\langle 93\rangle\)}}
\put(-2.7,-10){\makebox(0,0)[lc]{\(2;7\)}}

\put(-0.7,-10){\makebox(0,0)[lc]{\(2;2\), periodic root}}
\put(-0.7,-11.5){\makebox(0,0)[lc]{\(1;2\)}}

\put(1.3,-11.5){\makebox(0,0)[lc]{\(1;1\)}}
\put(1.3,-13){\makebox(0,0)[lc]{\(1;2\)}}

\put(3.3,-13){\makebox(0,0)[lc]{\(1;1\)}}
\put(4.3,-13.7){\makebox(0,0)[lt]{etc.}}
\put(3.3,-14.5){\makebox(0,0)[lc]{\(1;2\)}}

% non-metabelian vertices
\put(-8.7,-7){\makebox(0,0)[lc]{\(\langle 622\rangle\)}}
\put(-9,-7.5){\makebox(0,0)[lc]{\(\mathrm{E}.8\)}}
\put(-9,-8){\makebox(0,0)[lc]{\((3,3)\)}}

\put(-6.7,-8.5){\makebox(0,0)[lc]{\(2;8\)}}
\put(-7,-9){\makebox(0,0)[lc]{\(\mathrm{D}.10\)}}
\put(-7,-9.5){\makebox(0,0)[lc]{\((9,3)\)}}

\put(-4.7,-10){\makebox(0,0)[lc]{\(2;2\)}}
\put(-5,-10.5){\makebox(0,0)[lc]{\(\mathrm{D}.10\)}}
\put(-5,-11){\makebox(0,0)[lc]{\((27,3)\)}}

\put(-2.7,-11.5){\makebox(0,0)[lc]{\(1;1\)}}
\put(-3,-12){\makebox(0,0)[lc]{\(\mathrm{D}.10\)}}
\put(-3,-12.5){\makebox(0,0)[lc]{\((81,3)\)}}

\put(-0.7,-13){\makebox(0,0)[lc]{\(1;1\)}}
\put(-1,-13.5){\makebox(0,0)[lc]{\(\mathrm{D}.10\)}}
\put(-1,-14){\makebox(0,0)[lc]{\((243,3)\)}}

\put(1.3,-14.5){\makebox(0,0)[lc]{\(1;1\)}}
\put(1,-15){\makebox(0,0)[lc]{\(\mathrm{D}.10\)}}
\put(1,-15.5){\makebox(0,0)[lc]{\((729,3)\)}}

\put(3.3,-16){\makebox(0,0)[lc]{\(1;1\)}}
\put(3,-16.5){\makebox(0,0)[lc]{\(\mathrm{D}.10\)}}
\put(3,-17){\makebox(0,0)[lc]{\((2187,3)\)}}

%***********************************************
% legend
\put(-3,-0.6){\makebox(0,0)[lc]{Legend:}}

\put(-1.6,-0.7){\framebox(0.2,0.2){}}
\put(-1.5,-0.6){\circle*{0.1}}
\put(-1,-0.6){\makebox(0,0)[lc]{\(\ldots\) abelian}}

\put(-1.5,-1){\circle{0.2}}
\put(-1,-1){\makebox(0,0)[lc]{\(\ldots\) metabelian}}

\put(-1.5,-1.4){\circle*{0.2}}
\put(-1,-1.4){\makebox(0,0)[lc]{\(\ldots\) metabelian with bifurcation}}

\put(-1.6,-1.9){\framebox(0.2,0.2){}}
\put(-1,-1.8){\makebox(0,0)[lc]{\(\ldots\) non-metabelian Schur \(\sigma\)}}

\end{picture}

}

\end{figure}
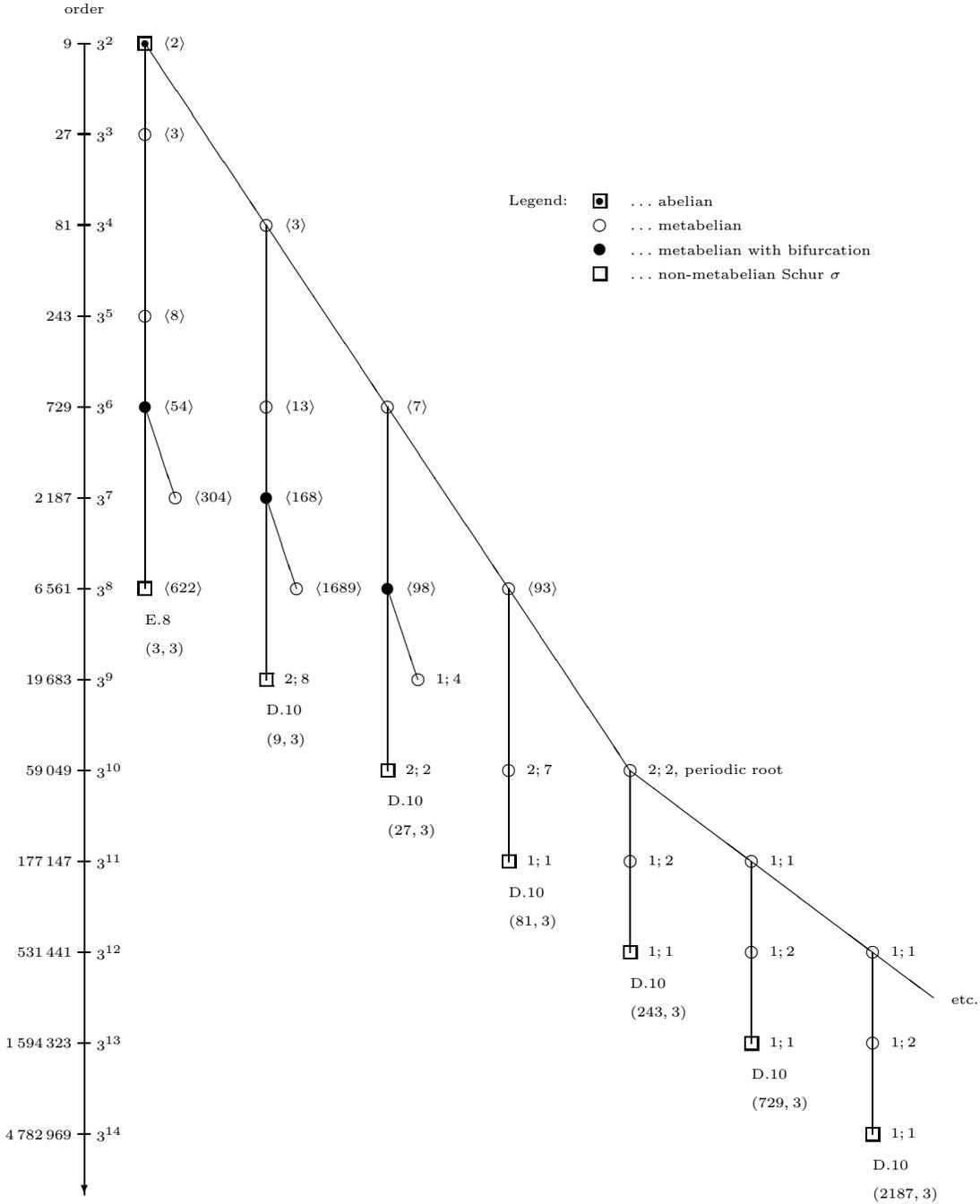

In the main theorem
\cite[\S\ 9, Thm. 12]{Ma2021}
of the previous work,
we provided evidence of another periodicity
of pairs of \textit{non-metabelian} Schur \(\sigma\)-groups \(G\)
with \(G/G^\prime\simeq (3^e,3)\), \(e\ge 5\),
and four types
\(\mathrm{D}.5\), \(\varkappa\sim (211;3)\),
\(\mathrm{C}.4\), \(\varkappa\sim (311;3)\),
\(\mathrm{D}.10\), \(\varkappa\sim (411;3)\), and
\(\mathrm{D}.6\), \(\varkappa\sim (123;1)\),
which is illustrated for one member of the pair of type \(\mathrm{D}.10\) by Figure
\ref{fig:SchurSigma5}.

%\newpage
%--------------------------------------------------------------------------------

%\noindent
Figure
\ref{fig:SchurSigma5}
shows that the construction process for the eight
non-metabelian Schur \(\sigma\)-groups \(G\) with order \(\#G=3^{7+e}\) and
punctured transfer kernel types \(\mathrm{D}.10\), \(\mathrm{C}.4\), \(\mathrm{D}.5\), and \(\mathrm{D}.6\), 
becomes increasingly difficult for the commutator quotients \(G/G^\prime\simeq (27,3)\), \((81,3)\), \((243,3)\).
For the commutator quotient \(G/G^\prime\simeq (729,3)\), however,
an \textbf{unexpected tranquilization} occurs,
and the construction process becomes settled with a \textbf{simple step size one periodicity}.

%\newpage
%--------------------------------------------------------------------------------

\begin{figure}[hb]
\caption{Schur \(\sigma\)-groups \(G\) with \(\varrho(G)\sim (3,3,3;3)\), \(G/G^\prime\simeq (3^e,3)\), \(2\le e\le 4\)}
\label{fig:SchurSigma}

{\tiny

\setlength{\unitlength}{0.9cm}
\begin{picture}(14,22.5)(-11,-19.5)

% scale of orders
\put(-11,2.5){\makebox(0,0)[cb]{order}}

\put(-11,2){\line(0,-1){21}}
\multiput(-11.1,2)(0,-1){22}{\line(1,0){0.2}}

%***********************************************

\put(-10.8,2){\makebox(0,0)[lc]{\(3^2\)}}
\put(-10.8,1){\makebox(0,0)[lc]{\(3^3\)}}
\put(-10.8,0){\makebox(0,0)[lc]{\(3^4\)}}
\put(-10.8,-1){\makebox(0,0)[lc]{\(3^5\)}}
\put(-10.8,-2){\makebox(0,0)[lc]{\(3^6\)}}
\put(-10.8,-3){\makebox(0,0)[lc]{\(3^7\)}}
\put(-10.8,-4){\makebox(0,0)[lc]{\(3^8\)}}
\put(-10.8,-5){\makebox(0,0)[lc]{\(3^9\)}}
\put(-10.8,-6){\makebox(0,0)[lc]{\(3^{10}\)}}
\put(-10.8,-7){\makebox(0,0)[lc]{\(3^{11}\)}}
\put(-10.8,-8){\makebox(0,0)[lc]{\(3^{12}\)}}
\put(-10.8,-9){\makebox(0,0)[lc]{\(3^{13}\)}}
\put(-10.8,-10){\makebox(0,0)[lc]{\(3^{14}\)}}
\put(-10.8,-11){\makebox(0,0)[lc]{\(3^{15}\)}}
\put(-10.8,-12){\makebox(0,0)[lc]{\(3^{16}\)}}
\put(-10.8,-13){\makebox(0,0)[lc]{\(3^{17}\)}}
\put(-10.8,-14){\makebox(0,0)[lc]{\(3^{18}\)}}
\put(-10.8,-15){\makebox(0,0)[lc]{\(3^{19}\)}}
\put(-10.8,-16){\makebox(0,0)[lc]{\(3^{20}\)}}
\put(-10.8,-17){\makebox(0,0)[lc]{\(3^{21}\)}}
\put(-10.8,-18){\makebox(0,0)[lc]{\(3^{22}\)}}
\put(-10.8,-19){\makebox(0,0)[lc]{\(3^{23}\)}}
%\put(-10.8,-20){\makebox(0,0)[lc]{\(3^{24}\)}}

%***********************************************
% vertices

\put(-11,-19){\vector(0,-1){1}}

% abelian root vertex
\put(-9.1,1.9){\framebox(0.2,0.2){}}
\put(-9,2){\circle*{0.1}}

% metabelian vertices
\put(-9,1){\circle{0.2}}
\put(-9,-1){\circle{0.2}}
\put(-9,-3){\circle*{0.2}}
\put(-8.5,-5){\circle{0.2}}

\put(-7,0){\circle{0.2}}
\put(-7,-2){\circle{0.2}}
\put(-7,-4){\circle*{0.2}}
\put(-6.5,-6){\circle{0.2}}

\put(-5,-3){\circle{0.2}}
\put(-5,-5){\circle*{0.2}}
\put(-4.5,-7){\circle{0.2}}

\put(-3,-6){\circle*{0.2}}
\put(-2.5,-8){\circle{0.2}}

% non-metabelian vertices
\put(-9.1,-7.1){\framebox(0.2,0.2){}}
\put(-9.1,-9.1){\framebox(0.2,0.2){}}
\put(-9.1,-13.1){\framebox(0.2,0.2){}}
\put(-9.1,-14.1){\framebox(0.2,0.2){}}
\put(-9.1,-16.1){\framebox(0.2,0.2){}}
\put(-9,-16){\circle{0.1}}

\put(-7.1,-8.1){\framebox(0.2,0.2){}}
\put(-7.1,-10.1){\framebox(0.2,0.2){}}
\put(-7.1,-14.1){\framebox(0.2,0.2){}}
\put(-7.1,-15.1){\framebox(0.2,0.2){}}
\put(-7.1,-17.1){\framebox(0.2,0.2){}}
\put(-7,-17){\circle{0.1}}

\put(-5.1,-9.1){\framebox(0.2,0.2){}}
\put(-5.1,-11.1){\framebox(0.2,0.2){}}
\put(-5.1,-15.1){\framebox(0.2,0.2){}}
\put(-5.1,-16.1){\framebox(0.2,0.2){}}
\put(-5.1,-18.1){\framebox(0.2,0.2){}}
\put(-5,-18){\circle{0.1}}

\put(-3.1,-10.1){\framebox(0.2,0.2){}}
\put(-3.1,-12.1){\framebox(0.2,0.2){}}
\put(-3.1,-16.1){\framebox(0.2,0.2){}}
\put(-3.1,-17.1){\framebox(0.2,0.2){}}
\put(-3.1,-19.1){\framebox(0.2,0.2){}}
\put(-3,-19){\circle{0.1}}

%\put(-1.1,-20.1){\framebox(0.2,0.2){}}
%\put(-1,-20){\circle{0.1}}

%***********************************************
% directed edges 

\put(-9,2){\line(0,-1){1}}
\put(-9,1){\line(0,-1){2}}
\put(-9,-1){\line(0,-1){2}}
\put(-9,-3){\line(0,-1){4}}
\put(-9,-3){\line(1,-4){0.5}}
\put(-9,-7){\line(0,-1){2}}
\put(-9,-9){\line(0,-1){4}}
\put(-9,-13){\line(0,-1){1}}
\put(-9,-14){\line(0,-1){2}}

\put(-9,2){\line(1,-1){2}}
\put(-7,0){\line(0,-1){2}}
\put(-7,-2){\line(0,-1){2}}
\put(-7,-4){\line(0,-1){4}}
\put(-7,-4){\line(1,-4){0.5}}
\put(-7,-8){\line(0,-1){2}}
\put(-7,-10){\line(0,-1){4}}
\put(-7,-14){\line(0,-1){1}}
\put(-7,-15){\line(0,-1){2}}

\put(-7,0){\line(2,-3){2}}
\put(-5,-3){\line(0,-1){2}}
\put(-5,-5){\line(0,-1){4}}
\put(-5,-5){\line(1,-4){0.5}}
\put(-5,-9){\line(0,-1){2}}
\put(-5,-11){\line(0,-1){4}}
\put(-5,-15){\line(0,-1){1}}
\put(-5,-16){\line(0,-1){2}}

\put(-5,-3){\line(2,-3){2}}
\put(-3,-6){\line(0,-1){4}}
\put(-3,-6){\line(1,-4){0.5}}
\put(-3,-10){\line(0,-1){2}}
\put(-3,-12){\line(0,-1){4}}
\put(-3,-16){\line(0,-1){1}}
\put(-3,-17){\line(0,-1){2}}

\put(-3,-6){\line(1,-1){2}}
\put(-0.8,-8){\makebox(0,0)[lc]{\((243,3)\), continued in Fig. \ref{fig:SchurSigma59}}}

%***********************************************
% identifiers in the SmallGroups Database

% abelian root vertex
\put(-8.8,2){\makebox(0,0)[lc]{\(\langle 2\rangle\)}}

% metabelian vertices
\put(-8.8,1){\makebox(0,0)[lc]{\(\langle 3\rangle\)}}
\put(-8.8,-1){\makebox(0,0)[lc]{\(\langle 3\rangle\)}}
\put(-8.8,-3){\makebox(0,0)[lc]{\(\langle 64\rangle\)}}
\put(-8.3,-5){\makebox(0,0)[lc]{\(2;36\)}}

\put(-6.8,0){\makebox(0,0)[lc]{\(\langle 3\rangle\)}}
\put(-6.8,-2){\makebox(0,0)[lc]{\(\langle 10\rangle\)}}
\put(-6.8,-4){\makebox(0,0)[lc]{\(\langle 165\rangle\)}}
\put(-6.3,-6){\makebox(0,0)[lc]{\(2;85\)}}

\put(-4.8,-3){\makebox(0,0)[lc]{\(\langle 3\rangle\)}}
\put(-4.8,-5){\makebox(0,0)[lc]{\(2;10\)}}
\put(-4.3,-7){\makebox(0,0)[lc]{\(2;88\)}}

\put(-2.8,-6){\makebox(0,0)[lc]{\(3;2\), the last bifurcation, \(\mathrm{cl}=4=e\)}}
\put(-2.3,-8){\makebox(0,0)[lc]{\(2;100\)}}

% non-metabelian vertices
\put(-8.8,-7){\makebox(0,0)[lc]{\(4;144\)}}
\put(-8.8,-9){\makebox(0,0)[lc]{\(2;41\)}}
\put(-8.8,-13){\makebox(0,0)[lc]{\(4;1\)}}
\put(-8.8,-14){\makebox(0,0)[lc]{\(1;1\)}}
\put(-8.8,-16){\makebox(0,0)[lc]{\(2;1\)}}
\put(-9,-16.5){\makebox(0,0)[lc]{\(\mathrm{F}.11\)}}
\put(-9,-17){\makebox(0,0)[lc]{\((3,3)\)}}

\put(-6.8,-8){\makebox(0,0)[lc]{\(4;29\)}}
\put(-6.8,-10){\makebox(0,0)[lc]{\(2;41\)}}
\put(-6.8,-14){\makebox(0,0)[lc]{\(4;1\)}}
\put(-6.8,-15){\makebox(0,0)[lc]{\(1;2\)}}
\put(-6.8,-17){\makebox(0,0)[lc]{\(2;1\)}}
\put(-7,-17.5){\makebox(0,0)[lc]{\(\mathrm{B}.18\)}}
\put(-7,-18){\makebox(0,0)[lc]{\((9,3)\)}}

\put(-4.8,-9){\makebox(0,0)[lc]{\(4;51\)}}
\put(-4.8,-11){\makebox(0,0)[lc]{\(2;32\)}}
\put(-4.8,-15){\makebox(0,0)[lc]{\(4;1\)}}
\put(-4.8,-16){\makebox(0,0)[lc]{\(1;1\)}}
\put(-4.8,-18){\makebox(0,0)[lc]{\(2;1\)}}
\put(-5,-18.5){\makebox(0,0)[lc]{\(\mathrm{B}.18\)}}
\put(-5,-19){\makebox(0,0)[lc]{\((27,3)\)}}

\put(-2.8,-10){\makebox(0,0)[lc]{\(4;81\)}}
\put(-2.8,-12){\makebox(0,0)[lc]{\(2;5\)}}
\put(-2.8,-16){\makebox(0,0)[lc]{\(4;1\)}}
\put(-2.8,-17){\makebox(0,0)[lc]{\(1;1\)}}
\put(-2.8,-19){\makebox(0,0)[lc]{\(2;1\)}}
\put(-3,-19.5){\makebox(0,0)[lc]{\(\mathrm{B}.18\)}}
\put(-3,-20){\makebox(0,0)[lc]{\((81,3)\)}}

%\put(-0.8,-20){\makebox(0,0)[lc]{\(2;1\) ???}}
%\put(-1,-20.5){\makebox(0,0)[lc]{\(\mathrm{B}.18\)}}
%\put(-1,-21){\makebox(0,0)[lc]{\((243,3)\)}}

%***********************************************
% legend
\put(-4,-0.6){\makebox(0,0)[lc]{Legend:}}

\put(-2.6,-0.7){\framebox(0.2,0.2){}}
\put(-2.5,-0.6){\circle*{0.1}}
\put(-2,-0.6){\makebox(0,0)[lc]{\(\ldots\) abelian}}

\put(-2.5,-1){\circle{0.2}}
\put(-2,-1){\makebox(0,0)[lc]{\(\ldots\) metabelian}}

\put(-2.5,-1.4){\circle*{0.2}}
\put(-2,-1.4){\makebox(0,0)[lc]{\(\ldots\) metabelian with bifurcation}}

\put(-2.6,-1.9){\framebox(0.2,0.2){}}
\put(-2,-1.8){\makebox(0,0)[lc]{\(\ldots\) non-metabelian}}

\put(-2.6,-2.3){\framebox(0.2,0.2){}}
\put(-2.5,-2.2){\circle{0.1}}
\put(-2,-2.2){\makebox(0,0)[lc]{\(\ldots\) non-metabelian Schur \(\sigma\)}}

\end{picture}

}

\end{figure}

%\newpage
%--------------------------------------------------------------------------------

\noindent
The investigation of periodic Schur \(\sigma\)-groups \(G\)
with \textit{moderate} rank distribution
\(\varrho(G)\sim (2,2,2;3)\) or \(\varrho(G)\sim (2,2,3;3)\)
was completed in
\cite{Ma2021}.
Although we were conscious that the difficulties will increase significantly,
the tree diagram in Figure
\ref{fig:SchurSigma5}
inspired us to look at cases with \textit{elevated} rank distribution
on \(21\) August \(2021\).
In Figure
\ref{fig:SchurSigma},
we see how large Schur \(\sigma\)-groups \(G\)
with logarithmic order \(\mathrm{lo}(G)=19+e\)
and commutator quotient \(G/G^\prime\simeq (3^e,3)\), \(1\le e\le 4\),
can be constructed with the aid of the \(p\)-group generation algorithm
\cite{Nm1977,Ob1990},
which is implemented in the ANUPQ package
\cite{GNO2006}
of the computational algebra system Magma
\cite{MAGMA2021}.
In these four cases,
the exponent \(e\) is not bigger than
the nilpotency class \(\mathrm{cl}(F)=4\) of the metabelian fork \(F\)
with bifurcation to non-metabelian vertices
\[
G\buildrel s=2 \over\longrightarrow
\pi(G)\buildrel s=1 \over\longrightarrow
\pi^2(G)\buildrel s=4 \over\longrightarrow
\pi^3(G)\buildrel s=2 \over\longrightarrow
\pi^4(G)\buildrel s=4 \over\longrightarrow
\pi^5(G)=F.
\]
They form the \textit{extremal root path} of the Schur \(\sigma\)-group \(G\),
which is weighted by the maximal step sizes \(s=\nu\)
equal to the \textit{nuclear rank} of the parent.
In this region, parents and \(p\)-parents coincide.

Figure
\ref{fig:SchurSigma}
for \(2\le e\le 4\),
which is continued by Figure
\ref{fig:SchurSigma59}
for \(4\le e\le 13\),
documents the stagnating state of our research enterprise
on \(31\) August \(2021\),
due to group theoretic problems.
The initial cases were still in the region where parents and \(p\)-parents coincide,
\[
\text{for } e = 3: \quad
\langle 2187,3\rangle \buildrel s=2 \over\longleftarrow \langle 2187,3\rangle-\#2;10
\buildrel s=4 \over\longleftarrow \langle 2187,3\rangle-\#2;10-\#4;51 \longleftarrow \quad \text{ etc.}
\]
\[
\text{for } e = 4: \quad
\langle 2187,3\rangle-\#3;2 \buildrel s=4 \over\longleftarrow \langle 2187,3\rangle-\#3;2-\#4;81 \longleftarrow \quad \text{ etc.}
\]
However, the case \(e = 5\) was outside of our reach already.
We tried to look at the descendant \(\langle 2187,3\rangle-\#3;2-\#5;1\), which has \(G/G^\prime \simeq (243,3)\),
but we got too big AQI of first order, namely
\((622,511,511;421)\) instead of \((621,511,511;421)\).

%Theorems were only given for \S\S\
%\ref{s:31},
%\ref{s:Imag31}
%and
%\ref{s:Imag41}.
%Similar theorems exist for \S\
%\ref{s:Imag21}
%but they were already presented in
%\cite{Ma2020b}.
%Theorems for \S\S\
%\ref{s:Imag51}
%and
%\ref{s:Imag61}
%could not be given,
%because the groups could not be constructed.
%(Now we can construct the groups but irregular behavior still prohibits the statement of theorems.)

%Our group theoretic problems are illuminated by Figure 
%\ref{fig:SchurSigma}.
At the commutator quotient \((81,3)=(3^e,3)\) with \(e=4\)
the exponent \(e\) overtakes the nilpotency class of the bifurcation \(\mathrm{cl}(F)=4\).
It was not clear if the bifurcation will vanish for \((243,3)\),
but eventually it turned out that \(B = \langle 2187,3\rangle-\#3;2\)
is simultaneous bifurcation for all commutator quotients \((3^e,3)\) with \(e\ge 4\).
It can thus be called a \textit{bifurcation of infinite order}.
%At the moment this was a mystery,
%since we were unable to see how to proceed further at \(\langle 2187,3\rangle-\#3;2\).
%It was (almost) clear that Schur \(\sigma\)-groups of minimal order \(3^{24}\) must exist,
%but how can they be constructed?

\bigskip
After a lot of trial and error
we succeeded in the construction of the desired
Schur \(\sigma\)-groups G with
logarithmic order \(24\),
\(G/G^\prime \simeq (243,3)\), i.e. \(e = 5\),
type \(\mathrm{B}.18\), \(\varkappa \sim (144;4)\),
AQI \(\alpha_1 \sim (621,511,511;421)\), and
\(\mathrm{sl} = 3\).
The mystery was solved
on \(06\) September \(2021\)
in the following way,
which finally lead to Figure
\ref{fig:SchurSigma59}
on \(13\) September \(2021\).
Let \(B := \langle 2187,3\rangle-\#3;2\).

In a first step,
we looked for the metabelianization \(M = G/G^{\prime\prime}\),
and we got two unique solutions:
\(M = B-\#2;93-\#1;i\) with \(i\in\lbrace 2,3\rbrace\).
%(We were not sure about
%\(B-\#2;92-\#1;i\) with \(i=2,3\),
%but it seems that the \(\sigma\)-automorphism is not relator inverting for step size \(s=1\).)

In a second step, 
we sought the non-metabelian Schur \(\sigma\)-group G.
There are \(15\) possible starting points,
\(B-\#4;k\) with \(23\le k\le 37\),
but only \(k\in\lbrace 24,26,28,30,31,33,37\rbrace\)
leads to Schur \(\sigma\)-groups with minimal \(\mathrm{lo}(G)=19+e\).
(The other values of \(k\) lead to \(M = B-\#2;92-\#1;i\) with \(i\in\lbrace 2,3\rbrace\).)
Exemplarily we take \(k=37\) in Figure
\ref{fig:SchurSigma59}.
The \textbf{classical root path} with respect to the usual lower central
\textbf{becomes disconnected}.
The first non-metabelian vertex is irregular,
\(B-\#4;37-\#1;i\) with \(i\in\lbrace 2,3\rbrace\).
It is isolated, since it has nuclear rank zero,
and thus is useless for the construction.
The remaining four non-metabelian vertices are regularly connected,
beginning at
\(B-\#4;37-\#3;j\) with \(j\in\lbrace 73,114\rbrace\).
We take \(j=73\) in Figure
\ref{fig:SchurSigma59},
which therefore illustrates a particular instance of the main Theorem
\ref{thm:Elevated37}.
The structure of the relevant tree diagrams
for the other five main Theorems
\ref{thm:Elevated24} --
\ref{thm:Elevated33}
is the same as in Figure
\ref{fig:SchurSigma59}.

\bigskip
Concerning the bifurcations, we have the following information:

\begin{theorem}
\label{thm:Bifurcations}
The bifurcations possess nearly identical pc-presentations:
there are in fact only three bifurcations,
\(B = \langle 6561,165\rangle\) for \((9,3)\),
\(B = \langle 2187,3\rangle-\#2;10\) for \((27,3)\), and
the \textbf{bifurcation of infinite order}
\(B = \langle 2187,3\rangle-\#3;2\) for any \((3^e,3)\) with \(e \ge 4\).
Denote some crucial commutators by
\(s_2 = \lbrack y,x\rbrack\),
\(s_3 = \lbrack s_2,x\rbrack\), \(t_3 = \lbrack s_2,y\rbrack\),
\(s_4 = \lbrack s_3,x\rbrack\), \(t_4 = \lbrack t_3,y\rbrack\),
\(s_5 = \lbrack s_4,x\rbrack\), \(t_5 = \lbrack t_4,y\rbrack\).
Then the polycyclic pc-presentation is given by
\begin{equation}
\label{eqn:Bifurcations}
B = \langle x,y \mid x^{3^e}=1, y^3=s_3s_4^2, s_2^3=s_4t_4^2, \lbrack x^3,y\rbrack =s_4t_4\rangle
\end{equation}
with \(e=2\), respectively \(e=3\), respectively \(e=4\).
\end{theorem}

\newpage
%--------------------------------------------------------------------------------

\begin{figure}[ht]
\caption{Schur \(\sigma\)-groups \(G\) with \(\varrho(G)\sim (3,3,3;3)\), \(G/G^\prime\simeq (3^e,3)\), \(4\le e\le 13\)}
\label{fig:SchurSigma59}

{\tiny

\setlength{\unitlength}{0.9cm}
\begin{picture}(14,22.5)(-10,-19)

% scale of orders
\put(-11,2.5){\makebox(0,0)[cb]{order}}

\put(-11,2){\line(0,-1){22}}
\multiput(-11.1,2)(0,-1){23}{\line(1,0){0.2}}

%***********************************************

\put(-10.8,2){\makebox(0,0)[lc]{\(3^{10}\)}}
\put(-10.8,1){\makebox(0,0)[lc]{\(3^{11}\)}}
\put(-10.8,0){\makebox(0,0)[lc]{\(3^{12}\)}}
\put(-10.8,-1){\makebox(0,0)[lc]{\(3^{13}\)}}
\put(-10.8,-2){\makebox(0,0)[lc]{\(3^{14}\)}}
\put(-10.8,-3){\makebox(0,0)[lc]{\(3^{15}\)}}
\put(-10.8,-4){\makebox(0,0)[lc]{\(3^{16}\)}}
\put(-10.8,-5){\makebox(0,0)[lc]{\(3^{17}\)}}
\put(-10.8,-6){\makebox(0,0)[lc]{\(3^{18}\)}}
\put(-10.8,-7){\makebox(0,0)[lc]{\(3^{19}\)}}
\put(-10.8,-8){\makebox(0,0)[lc]{\(3^{20}\)}}
\put(-10.8,-9){\makebox(0,0)[lc]{\(3^{21}\)}}
\put(-10.8,-10){\makebox(0,0)[lc]{\(3^{22}\)}}
\put(-10.8,-11){\makebox(0,0)[lc]{\(3^{23}\)}}
\put(-10.8,-12){\makebox(0,0)[lc]{\(3^{24}\)}}
\put(-10.8,-13){\makebox(0,0)[lc]{\(3^{25}\)}}
\put(-10.8,-14){\makebox(0,0)[lc]{\(3^{26}\)}}
\put(-10.8,-15){\makebox(0,0)[lc]{\(3^{27}\)}}
\put(-10.8,-16){\makebox(0,0)[lc]{\(3^{28}\)}}
\put(-10.8,-17){\makebox(0,0)[lc]{\(3^{29}\)}}
\put(-10.8,-18){\makebox(0,0)[lc]{\(3^{30}\)}}
\put(-10.8,-19){\makebox(0,0)[lc]{\(3^{31}\)}}
\put(-10.8,-20){\makebox(0,0)[lc]{\(3^{32}\)}}

%***********************************************
% vertices

\put(-11,-20){\vector(0,-1){1}}

% metabelian vertices

\put(-10,2){\circle*{0.2}}
\put(-9.5,0){\circle{0.2}}
\put(-8,0){\circle{0.2}}
\put(-8,-1){\circle{0.2}}
\put(-7,-1){\circle{0.2}}
\put(-7,-2){\circle{0.2}}

% non-metabelian vertices

\put(-10.1,-2.1){\framebox(0.2,0.2){}}
\put(-10.1,-4.1){\framebox(0.2,0.2){}}
\put(-10.1,-8.1){\framebox(0.2,0.2){}}
\put(-10.1,-9.1){\framebox(0.2,0.2){}}
\put(-10.1,-11.1){\framebox(0.2,0.2){}}
\put(-10,-11){\circle{0.1}}

\put(-5.1,-2.1){\framebox(0.2,0.2){}}
\put(-5.1,-5.1){\framebox(0.2,0.2){}}
\put(-5.1,-9.1){\framebox(0.2,0.2){}}
\put(-5.1,-10.1){\framebox(0.2,0.2){}}
\put(-5.1,-12.1){\framebox(0.2,0.2){}}
\put(-5,-12){\circle{0.1}}

\put(-4.1,-5.1){\framebox(0.2,0.2){}}
\put(-4.1,-9.1){\framebox(0.2,0.2){}}
\put(-4.1,-11.1){\framebox(0.2,0.2){}}
\put(-4.1,-13.1){\framebox(0.2,0.2){}}
\put(-4,-13){\circle{0.1}}

\put(-3.1,-9.1){\framebox(0.2,0.2){}}
\put(-3.1,-11.1){\framebox(0.2,0.2){}}
\put(-3.1,-13.1){\framebox(0.2,0.2){}}
\put(-3.1,-14.1){\framebox(0.2,0.2){}}
\put(-3,-14){\circle{0.1}}

\put(-2.1,-11.1){\framebox(0.2,0.2){}}
\put(-2.1,-13.1){\framebox(0.2,0.2){}}
\put(-2.1,-14.1){\framebox(0.2,0.2){}}
\put(-2.1,-15.1){\framebox(0.2,0.2){}}
\put(-2,-15){\circle{0.1}}

\put(-1.1,-13.1){\framebox(0.2,0.2){}}
\put(-1.1,-14.1){\framebox(0.2,0.2){}}
\put(-1.1,-15.1){\framebox(0.2,0.2){}}
\put(-1.1,-16.1){\framebox(0.2,0.2){}}
\put(-1,-16){\circle{0.1}}

\put(-0.1,-14.1){\framebox(0.2,0.2){}}
\put(-0.1,-15.1){\framebox(0.2,0.2){}}
\put(-0.1,-16.1){\framebox(0.2,0.2){}}
\put(-0.1,-17.1){\framebox(0.2,0.2){}}
\put(0,-17){\circle{0.1}}

\put(0.9,-15.1){\framebox(0.2,0.2){}}
\put(0.9,-16.1){\framebox(0.2,0.2){}}
\put(0.9,-17.1){\framebox(0.2,0.2){}}
\put(0.9,-18.1){\framebox(0.2,0.2){}}
\put(1,-18){\circle{0.1}}

\put(1.9,-16.1){\framebox(0.2,0.2){}}
\put(1.9,-17.1){\framebox(0.2,0.2){}}
\put(1.9,-18.1){\framebox(0.2,0.2){}}
\put(1.9,-19.1){\framebox(0.2,0.2){}}
\put(2,-19){\circle{0.1}}

\put(2.9,-17.1){\framebox(0.2,0.2){}}
\put(2.9,-18.1){\framebox(0.2,0.2){}}
\put(2.9,-19.1){\framebox(0.2,0.2){}}
\put(2.9,-20.1){\framebox(0.2,0.2){}}
\put(3,-20){\circle{0.1}}

%***********************************************
% directed edges 

\put(-10,2){\line(0,-1){4}}
\put(-10,2){\line(1,-4){0.5}}
\put(-10,-2){\line(0,-1){2}}
\put(-10,-4){\line(0,-1){4}}
\put(-10,-8){\line(0,-1){1}}
\put(-10,-9){\line(0,-1){2}}

\put(-10,2){\line(1,-1){2}}
\put(-8,0){\line(0,-1){1}}
\put(-8,0){\line(1,-1){1}}
\put(-7,-1){\line(0,-1){1}}

\put(-10,2){\line(5,-4){5}}
\put(-5,-2){\line(0,-1){3}}
\put(-5,-5){\line(0,-1){4}}
\put(-5,-9){\line(0,-1){1}}
\put(-5,-10){\line(0,-1){2}}

\put(-5,-2){\line(1,-3){1}}
\put(-4,-5){\line(0,-1){4}}
\put(-4,-9){\line(0,-1){2}}
\put(-4,-11){\line(0,-1){2}}

\put(-4,-5){\line(1,-4){1}}
\put(-3,-9){\line(0,-1){2}}
\put(-3,-11){\line(0,-1){2}}
\put(-3,-13){\line(0,-1){1}}

\put(-3,-9){\line(1,-2){1}}
\put(-2,-11){\line(0,-1){2}}
\put(-2,-13){\line(0,-1){1}}
\put(-2,-14){\line(0,-1){1}}

\put(-2,-11){\line(1,-2){1}}
\put(-1,-13){\line(0,-1){1}}
\put(-1,-14){\line(0,-1){1}}
\put(-1,-15){\line(0,-1){1}}

\put(-1,-13){\line(1,-1){1}}
\put(0,-14){\line(0,-1){1}}
\put(0,-15){\line(0,-1){1}}
\put(0,-16){\line(0,-1){1}}

\put(0,-14){\line(1,-1){1}}
\put(1,-15){\line(0,-1){1}}
\put(1,-16){\line(0,-1){1}}
\put(1,-17){\line(0,-1){1}}

\put(1,-15){\line(1,-1){1}}
\put(2,-16){\line(0,-1){1}}
\put(2,-17){\line(0,-1){1}}
\put(2,-18){\line(0,-1){1}}

\put(2,-16){\line(1,-1){1}}
\put(3,-17){\line(0,-1){1}}
\put(3,-18){\line(0,-1){1}}
\put(3,-19){\line(0,-1){1}}

%***********************************************
% identifiers in the ANUPQ package

% metabelian vertices

\put(-9.8,2){\makebox(0,0)[lc]{\(3;2\), the last bifurcation, \(\mathrm{cl}=4=e\), continued from Fig. \ref{fig:SchurSigma}}}
\put(-9.3,0){\makebox(0,0)[lc]{\(2;100\)}}
\put(-9.5,-0.3){\makebox(0,0)[lc]{\((81,3)\)}}
\put(-7.8,0){\makebox(0,0)[lc]{\(2;93\), metabelian periodic root}}
\put(-7.8,-1){\makebox(0,0)[lc]{\(1;2\)}}
\put(-8,-1.3){\makebox(0,0)[lc]{\((243,3)\)}}
\put(-6.8,-1){\makebox(0,0)[lc]{\(1;1\)}}
\put(-6.8,-2){\makebox(0,0)[lc]{\(1;2\)}}
\put(-7,-2.3){\makebox(0,0)[lc]{\((729,3)\)}}
\put(-6.2,-1.7){\makebox(0,0)[lc]{etc.}}

% non-metabelian vertices

\put(-9.8,-2){\makebox(0,0)[lc]{\(4;81\)}}
\put(-9.8,-4){\makebox(0,0)[lc]{\(2;5\)}}
\put(-9.8,-8){\makebox(0,0)[lc]{\(4;1\)}}
\put(-9.8,-9){\makebox(0,0)[lc]{\(1;1\)}}
\put(-9.8,-11){\makebox(0,0)[lc]{\(2;1\)}}
\put(-10,-11.5){\makebox(0,0)[lc]{\(\mathrm{B}.18\)}}
\put(-10,-12){\makebox(0,0)[lc]{\((81,3)\)}}

\put(-4.8,-2){\makebox(0,0)[lc]{\(4;37\), situation in Theorem \ref{thm:Elevated37}}}
\put(-4.8,-5){\makebox(0,0)[lc]{\(3;73\)}}
\put(-4.8,-9){\makebox(0,0)[lc]{\(4;1\)}}
\put(-4.8,-10){\makebox(0,0)[lc]{\(1;2\)}}
\put(-4.8,-12){\makebox(0,0)[lc]{\(2;1\)}}
\put(-5,-12.5){\makebox(0,0)[rc]{\(\mathrm{B}.18\)}}
\put(-5,-13){\makebox(0,0)[rc]{\((243,3)\)}}

\put(-3.8,-5){\makebox(0,0)[lc]{\(3;32\)}}
\put(-3.8,-9){\makebox(0,0)[lc]{\(4;10\)}}
\put(-3.8,-11){\makebox(0,0)[lc]{\(2;2\)}}
\put(-3.8,-13){\makebox(0,0)[lc]{\(2;1\)}}
\put(-4,-13.5){\makebox(0,0)[rc]{\(\mathrm{B}.18\)}}
\put(-4,-14){\makebox(0,0)[rc]{\((729,3)\)}}

\put(-2.8,-9){\makebox(0,0)[lc]{\(4;a\) with \(a=1\)}}
\put(-2.8,-11){\makebox(0,0)[lc]{\(2;4\)}}
\put(-2.8,-13){\makebox(0,0)[lc]{\(2;1\)}}
\put(-2.8,-14){\makebox(0,0)[lc]{\(1;1\)}}
\put(-3,-14.5){\makebox(0,0)[rc]{\(\mathrm{B}.18\)}}
\put(-3,-15){\makebox(0,0)[rc]{\((2187,3)\)}}

\put(-1.8,-11){\makebox(0,0)[lc]{\(2;\tilde{a}\) with \(\tilde{a}=2\)}}
\put(-1.8,-13){\makebox(0,0)[lc]{\(2;4\)}}
\put(-1.8,-14){\makebox(0,0)[lc]{\(1;1\)}}
\put(-1.8,-15){\makebox(0,0)[lc]{\(1;1\)}}
\put(-2,-15.5){\makebox(0,0)[rc]{\(\mathrm{B}.18\)}}
\put(-2,-16){\makebox(0,0)[rc]{\((6561,3)\)}}

\put(-0.8,-13){\makebox(0,0)[lc]{\(2;b\) with \(b=1\), non-metabelian periodic root \(W_{1,1}\)}}
\put(-0.8,-14){\makebox(0,0)[lc]{\(1;2\)}}
\put(-0.8,-15){\makebox(0,0)[lc]{\(1;1\)}}
\put(-0.8,-16){\makebox(0,0)[lc]{\(1;1\)}}
\put(-1,-16.5){\makebox(0,0)[rc]{\(\mathrm{B}.18\)}}
\put(-1,-17){\makebox(0,0)[rc]{\((19683,3)\)}}

\put(0.2,-14){\makebox(0,0)[lc]{\(1;1\)}}
\put(0.2,-15){\makebox(0,0)[lc]{\(1;2\)}}
\put(0.2,-16){\makebox(0,0)[lc]{\(1;1\)}}
\put(0.2,-17){\makebox(0,0)[lc]{\(1;1\)}}
\put(0,-17.5){\makebox(0,0)[rc]{\(\mathrm{B}.18\)}}
\put(0,-18){\makebox(0,0)[rc]{\((59049,3)\)}}

\put(1.2,-15){\makebox(0,0)[lc]{\(1;1\)}}
\put(1.2,-16){\makebox(0,0)[lc]{\(1;2\)}}
\put(1.2,-17){\makebox(0,0)[lc]{\(1;1\)}}
\put(1.2,-18){\makebox(0,0)[lc]{\(1;1\)}}
\put(1,-18.5){\makebox(0,0)[rc]{\(\mathrm{B}.18\)}}
\put(1,-19){\makebox(0,0)[rc]{\((177147,3)\)}}

\put(2.2,-16){\makebox(0,0)[lc]{\(1;1\)}}
\put(2.2,-17){\makebox(0,0)[lc]{\(1;2\)}}
\put(2.2,-18){\makebox(0,0)[lc]{\(1;1\)}}
\put(2.2,-19){\makebox(0,0)[lc]{\(1;1\)}}
\put(2,-19.5){\makebox(0,0)[rc]{\(\mathrm{B}.18\)}}
\put(2,-20){\makebox(0,0)[rc]{\((531441,3)\)}}

\put(3.2,-17){\makebox(0,0)[lc]{\(1;1\)}}
\put(3.2,-18){\makebox(0,0)[lc]{\(1;2\)}}
\put(3.2,-19){\makebox(0,0)[lc]{\(1;1\)}}
\put(3.2,-20){\makebox(0,0)[lc]{\(1;1\)}}
\put(3,-20.5){\makebox(0,0)[rc]{\(\mathrm{B}.18\)}}
\put(3,-21){\makebox(0,0)[rc]{\((1594323,3)\)}}

\put(4.2,-17.7){\makebox(0,0)[lc]{etc.}}

%***********************************************
% legend

\put(-4,1.4){\makebox(0,0)[lc]{Legend:}}

\put(-2.5,1.4){\circle{0.2}}
\put(-2,1.4){\makebox(0,0)[lc]{\(\ldots\) metabelian}}

\put(-2.5,1){\circle*{0.2}}
\put(-2,1){\makebox(0,0)[lc]{\(\ldots\) metabelian with bifurcation}}

\put(-2.6,0.5){\framebox(0.2,0.2){}}
\put(-2,0.6){\makebox(0,0)[lc]{\(\ldots\) non-metabelian}}

\put(-2.6,0.1){\framebox(0.2,0.2){}}
\put(-2.5,0.2){\circle{0.1}}
\put(-2,0.2){\makebox(0,0)[lc]{\(\ldots\) non-metabelian Schur \(\sigma\)}}

\end{picture}

}

\end{figure}

\newpage
%--------------------------------------------------------------------------------

\section{Imaginary quadratic fields \(K\) with \(\mathrm{Cl}_3(K)\simeq C_{243}\times C_3\)}
\label{s:Imag51}

\noindent
The \(1784\) imaginary quadratic fields \(K=\mathbb{Q}(\sqrt{d})\)
with fundamental discriminants \(-60\,000\,000<d<0\) and
\(3\)-class group \(\mathrm{Cl}_3(K)\simeq C_{243}\times C_3\)
were computed by means of the computational algebra system Magma
\cite{MAGMA2021}
in \(451\,227\) seconds of CPU time, that is nearly a full week.
In Table
\ref{tbl:Type51},
the first seven cases with punctured capitulation type \(\mathrm{B}.18\),
\(\varkappa(K)\sim (144;4)\),
are listed.
The abelian quotient invariants \(\alpha_1(K)\) of first order
of only six of them are \textit{uni-polarized} and in the \textit{ground state}.

%\newpage
%--------------------------------------------------------------------------------

\renewcommand{\arraystretch}{1.1}
\begin{table}[ht]
\caption{Seven fields \(K=\mathbb{Q}(\sqrt{d})\) with \(\mathrm{Cl}_3(K)\simeq C_{243}\times C_3\) and \(\varkappa(K)\sim (144;4)\)}
\label{tbl:Type51}
\begin{center}
\begin{tabular}{|r|r|c||c|c|}
\hline
 No.     & \(d\)             & factors           & \(\alpha_1(K)\)       & remark \\
\hline
  \(60\) &  \(-5\,629\,151\) & \(11,631,811\)    & \((621,511,511;421)\) & \\
  \(65\) &  \(-5\,702\,003\) & prime             & \((621,511,511;421)\) & \\
  \(73\) &  \(-6\,124\,411\) & prime             & \((621,511,511;421)\) & \\
  \(77\) &  \(-6\,219\,188\) & \(2,1\,554\,797\) & \((621,511,511;421)\) & \\
 \(116\) &  \(-8\,513\,951\) & prime             & \((621,511,511;432)\) & bi-polarized \\
 \(149\) & \(-10\,401\,044\) & \(2,41,63\,421\)  & \((621,511,511;421)\) & \\
 \(155\) & \(-10\,607\,215\) & \(5,2\,121\,443\) & \((621,511,511;421)\) & \\
\hline
\end{tabular}
\end{center}
\end{table}

%\newpage
%--------------------------------------------------------------------------------

\noindent
In Table
\ref{tbl:SecondAQI51},
we give the abelian quotient invariants \(\alpha_2(K)\) of second order
of the six fields in the uni-polarized ground state
contained in Table
\ref{tbl:Type51}.
The general structure of \(\alpha_2(K)\) is the following
\begin{equation}
\label{eqn:SecondAQI51}
\alpha_2(K)=\lbrack 51;(621;52111,D_1),(511;52111,D_2),(511;52111,D_3);(421;52111,D_4)\rbrack
\end{equation}
where each dodecuplet \(D_i\) usually consists of a triplet and a nonet.
Only the constitution of \(D_4\) is occasionally irregular.

%\newpage
%--------------------------------------------------------------------------------

\renewcommand{\arraystretch}{1.1}
\begin{table}[ht]
\caption{Details for six fields \(K=\mathbb{Q}(\sqrt{d})\) in Table \ref{tbl:Type51}}
\label{tbl:SecondAQI51}
\begin{center}
\begin{tabular}{|r||c|c|c|c|c|}
\hline
 No.     & \(D_1\)            & \(D_2\)             & \(D_3\)             & \(D_4\)                        & remark \\
\hline
  \(60\) & \((6211)^3(62)^9\) & \((621)^3(62)^9\)   & \((5211)^3(521)^9\) & \((531)^3(431)^6(422)^2(332)\) & irregular \\
  \(65\) & \(\) & \(\) & \(\) & \(\) & Magma int. err. \\
  \(73\) & \((6211)^3(62)^9\) & \((6111)^3(62)^9\) & \(\) & \((531)^3(521)^8(422)\) & Magma int. err. \\
  \(77\) & \(\)               & \((621)^3(62)^9\)  & \(\) & \((531)^3(521)^8(422)\) & Magma int. err. \\
 \(149\) & \(\) & \(\) & \(\) & \(\) & Magma int. err. \\
 \(155\) & \(\) & \(\) & \(\) & \(\) & Magma int. err. \\
\hline
\end{tabular}
\end{center}
\end{table}

%\newpage
%--------------------------------------------------------------------------------

\section{Imaginary quadratic fields \(K\) with \(\mathrm{Cl}_3(K)\simeq C_{729}\times C_3\)}
\label{s:Imag61}

\noindent
The \(263\) imaginary quadratic fields \(K=\mathbb{Q}(\sqrt{d})\)
with fundamental discriminants \(-60\,000\,000<d<0\) and
\(3\)-class group \(\mathrm{Cl}_3(K)\simeq C_{729}\times C_3\)
were computed by means of the computational algebra system Magma
\cite{MAGMA2021}
in \(411\,074\) seconds of CPU time, that is nearly a full week.
In Table
\ref{tbl:Type61},
the eleven cases with punctured capitulation type \(\mathrm{B}.18\),
\(\varkappa(K)\sim (144;4)\),
are listed.
The abelian quotient invariants \(\alpha_1(K)\) of first order
of only eight of them are \textit{uni-polarized} and in the \textit{ground state}.

%\newpage
%--------------------------------------------------------------------------------

\renewcommand{\arraystretch}{1.1}
\begin{table}[ht]
\caption{Eleven fields \(K=\mathbb{Q}(\sqrt{d})\) with \(\mathrm{Cl}_3(K)\simeq C_{729}\times C_3\) and \(\varkappa(K)\sim (144;4)\)}
\label{tbl:Type61}
\begin{center}
\begin{tabular}{|r|r|c||c|c|}
\hline
 No.     & \(d\)             & factors             & \(\alpha_1(K)\)       & remark \\
\hline
   \(9\) &  \(-8\,716\,319\) & \(2\,111,4\,129\)   & \((721,611,611;521)\) & \\
  \(17\) & \(-11\,598\,911\) & \(19,610\,469\)     & \((721,611,611;521)\) & \\
  \(28\) & \(-17\,054\,671\) & prime               & \((732,611,611;521)\) & first excited state \\
  \(94\) & \(-32\,670\,951\) & \(3,10\,890\,317\)  & \((721,611,611;521)\) & \\
 \(133\) & \(-38\,393\,396\) & \(2,9\,598\,349\)   & \((721,611,611;521)\) & \\
 \(141\) & \(-39\,551\,231\) & \(17,283,8\,221\)   & \((721,611,611;543)\) & highly bi-polarized \\
 \(144\) & \(-39\,948\,359\) & \(11,719,5\,051\)   & \((721,611,611;521)\) & \\
 \(197\) & \(-50\,631\,279\) & \(3,293,57\,601\)   & \((721,611,611;521)\) & \\
 \(198\) & \(-50\,963\,071\) & \(439,116\,089\)    & \((721,611,611;521)\) & \\
 \(242\) & \(-57\,507\,455\) & \(5,11\,501\,491\)  & \((721,611,611;521)\) & \\
 \(247\) & \(-58\,142\,996\) & \(2,14\,535\,749\)  & \((721,611,611;532)\) & bi-polarized \\
\hline
\end{tabular}
\end{center}
\end{table}

%\newpage
%--------------------------------------------------------------------------------

\noindent
In Table
\ref{tbl:SecondAQI61},
we give the abelian quotient invariants \(\alpha_2(K)\) of second order
of the eight fields in the uni-polarized ground state
contained in Table
\ref{tbl:Type61}.
The general structure of \(\alpha_2(K)\) is the following
\begin{equation}
\label{eqn:SecondAQI61}
\alpha_2(K)=\lbrack 61;(721;62111,D_1),(611;62111,D_2),(611;62111,D_3);(521;62111,D_4)\rbrack
\end{equation}
where each dodecuplet \(D_i\) usually consists of a triplet and a nonet.
Only the constitution of \(D_4\) is frequently (or even always) irregular.

%\newpage
%--------------------------------------------------------------------------------

\renewcommand{\arraystretch}{1.1}
\begin{table}[ht]
\caption{Details for eight fields \(K=\mathbb{Q}(\sqrt{d})\) in Table \ref{tbl:Type61}}
\label{tbl:SecondAQI61}
\begin{center}
\begin{tabular}{|r||c|c|c|c|c|}
\hline
 No.     & \(D_1\)            & \(D_2\)             & \(D_3\)            & \(D_4\)                         & remark \\
\hline
   \(9\) & \((7211)^3(72)^9\) & \((721)^3(621)^9\)  & \((721)^3(621)^9\) & \((5311)^3(531)^6(522)^2(432)\) & irregular \\
  \(17\) & \((7211)^3(72)^9\) & \((6211)^3(621)^9\) & \((721)^3(72)^9\)  & \((631)^3(531)^6(522)^2(432)\)  & irregular \\
  \(94\) & \(\) & \(\) & \(\) & \(\) & Magma int. err. \\
 \(133\) & \(\) & \(\) & \(\) & \(\) & Magma int. err. \\
 \(144\) & \(\) & \(\) & \(\) & \(\) & Magma int. err. \\
 \(197\) & \(\) & \(\) & \(\) & \((631)^3(531)^6(522)^2(432)\)  & irregular \\
 \(198\) & \(\) & \(\) & \(\) & \((5321)^3(531)^6(522)^2(432)\) & irregular \\
 \(242\) & \((7211)^3(72)^9\) & \((7111)^3(72)^9\)  & \(\)               & \((631)^3(531)^6(522)^2(432)\)  & irregular \\
\hline
\end{tabular}
\end{center}
\end{table}

%\newpage
%--------------------------------------------------------------------------------

\begin{remark}
\label{rmk:Type71And81}
In \S\
\ref{s:Layout}
we have mentioned that it is rather hopeless to continue the search for
imaginary quadratic fields \(K=\mathbb{Q}(\sqrt{d})\)
with bigger \(3\)-class groups \(\mathrm{Cl}_3(K)\simeq C_{3^e}\times C_3\) for \(e\ge 7\).
Firstly because of the immense amount of required CPU-time, and
secondly in view of Magma internal errors which occur with increasing frequency
during the computation of abelian type invariants \(\alpha_2(K)\) of the second order
for nonic relative extensions \(L/K\) with absolute degree \(18\).

Nevertheless, we mention some interesting observations in experiments with \(e=7\) and \(e=8\).
Concerning \(e=7\), we found six ground states of type \(\mathrm{B}.18\)
with \(\alpha_1(K)\sim (821,711,711;621)\)
for \(d\in\lbrace -37\,648\,463, -42\,705\,359, -122\,519\,927, -138\,616\,719, -154\,511\,167, -193\,538\,383\rbrace\),
and a bipolarization with \(\alpha_1(K)\sim (821,711,711;632)\)
for \(d=-206\,130\,371\). Five of these seven discriminants are prime.
The first two minimal hits of the desired \(3\)-class group are the primes \\
\(d=-32\,681\,951\) with \(\alpha_1(K)\sim (81,81,821;711)\) and type \(\mathrm{D}.5\), \(\varkappa(K)\sim (112;3)\), \\
\(d=-35\,574\,431\) with \(\alpha_1(K)\sim (81,81,711;711)\) and type \(\mathrm{D}.11\), \(\varkappa(K)\sim (124;1)\).

Concerning \(e=8\), we were at least able to discover three minimal hits
of the desired \(3\)-class group, though not of type \(\mathrm{B}.18\), \(\varkappa(K)\sim (144;4)\).
All discriminants are prime: \\
\(d=-98\,311\,919\) with \(\alpha_1(K)\sim (91,91,932;811)\) a first excited state of type \(\mathrm{D}.5\), \(\varkappa(K)\sim (112;3)\), \\
\(d=-201\,210\,239\) with \(\alpha_1(K)\sim (91,91,811;811)\) and type \(\mathrm{D}.11\), \(\varkappa(K)\sim (124;1)\), and \\
\(d=-209\,606\,759\) with \(\alpha_1(K)\sim (91,91,91;822)\) and type \(\mathrm{D}.6\), \(\varkappa(K)\sim (123;1)\).
\end{remark}

%\newpage
%--------------------------------------------------------------------------------

\section{A general theorem}
\label{s:General}

\noindent
The previous sections with concrete results for
various fixed values of the exponent \(2\le e\le 20\)
in the non-elementary bicyclic commutator quotient \(G/G^\prime\simeq (3^e,3)\hat{=}(e1)\)
suggest the following generalization with upper bound \(B:=20\).

\begin{theorem}
\label{thm:GeneralAQI}
In dependence on the exponent \(3\le e\le B\),
the abelian quotient invariants \(\alpha_2(G)\) of second order
of finite Schur \(\sigma\)-groups \(G\) with
commutator quotient \(G/G^\prime\simeq (e1)\),
punctured transfer kernel type \(\mathrm{B}.18\),
\(\varkappa(K)\sim (144;4)\),
and logarithmic order \(\mathrm{lo}(G)=19+e\)
are given by
\begin{equation}
\label{eqn:GeneralAQI}
\alpha_2(G)=\lbrack e1;((e+1)21;e2111,D_1),(e11;e2111,D_2),(e11;e2111,D_3);((e-1)21;e2111,D_4)\rbrack,
\end{equation}
where each dodecuplet \(D_i\) consists of a triplet \(T_i^3\) and a nonet \(N_i^9\),
except for \(i=4\), where the nonet \(N_4^9\) is replaced by an octet \(O_4^8\) and a singlet \(S_4\):
\begin{equation}
\label{eqn:Multiplets}
\begin{aligned}
T_1 &\in\lbrace (e+1)211,(e+1)1111\rbrace, \quad N_1=(e+1)2, \\
T_i &\in\lbrace (e+1)21,(e+1)111\rbrace, \quad N_i\in\lbrace (e+1)2,e21\rbrace, \quad \text{ for } 2\le i\le 3, \\
T_4 &\in\lbrace e31,e211\rbrace, \quad O_4=e21, \quad S_4=(e-1)22.
\end{aligned}
\end{equation}
\end{theorem}

\begin{conjecture}
\label{cnj:GeneralAQI}
Theorem
\ref{thm:GeneralAQI}
remains true for any upper bound \(B\ge 21\).
\end{conjecture}

\begin{remark}
\label{rmk:GeneralAQI}
\(T_i=e211\) for \(2\le i\le 3\) can also occur
but it leads to bigger logarithmic order \(\mathrm{lo}(G)>19+e\).
The same is true for a nonet \(N_4^9\) with \(N_4=e21\) in the fourth dodecuplet \(D_4\).
Theorem
\ref{thm:GeneralAQI}
was stated on \(31\) August \(2021\).
\end{remark}

%\newpage
%--------------------------------------------------------------------------------

\section{Conclusion}
\label{s:TheEnd}

\noindent
In our invited key note
\cite{Ma2020a}
at the \(3\)rd International Conference on Mathematics and its Applications (ICMA)
Casablanca, \(28\) February \(2020\),
we offered supervision of a Ph.D. thesis about
\(3\)-groups with non-elementary bicyclic commutator quotient
to the young researchers who listened to our talk and presentation
with vigilance.
That was eighteen months ago, immediately before
the breakout of the worldwide Corona crisis,
which prohibited any further scientific collaboration with personal contact.
In the present article and its predecessor
\cite{Ma2021}
we actually wrote this \lq\lq thesis\rq\rq\ ourselves,
thereby discovering several groundbreaking and totally unexpected simple periodicities.

%\newpage
%--------------------------------------------------------------------------------

\section{Acknowledgements}
\label{s:Gratifications}

\noindent
The author thanks Professor M. F. Newman
from the Australian National University in Canberra, Australian Capital Territory,
for his encouragement and aid during the preparation of this article.

The author acknowledges
financial support by the Austrian Science Fund (FWF): projects J0497-PHY, P26008-N25
and by the Research Executive Agency of the European Union (EUREA).

%\newpage
%--------------------------------------------------------------------------------

%--------------------------------------------------------------------------------


\begin{thebibliography}{XX}
%
\bibitem{Ag1998}
M. Arrigoni,
\textit{On Schur \(\sigma\)-groups},
Math. Nachr.
\textbf{192}
(1998),
71--89.
%
\bibitem{Ar1927}
E. Artin,
\textit{Beweis des allgemeinen Reziprozit\"atsgesetzes},
Abh. Math. Sem. Univ. Hamburg
\textbf{5}
(1927),
353--363.
%
\bibitem{Ar1929}
E. Artin,
\textit{Idealklassen in Oberk\"orpern und allgemeines Reziprozit\"atsgesetz},
Abh. Math. Sem. Univ. Hamburg
\textbf{7}
(1929),
46--51.
%
\bibitem{AHL1977}
J. A. Ascione, G. Havas, and C. R. Leedham-Green,
\textit{A computer aided classification of certain groups of prime power order},
Bull. Austral. Math. Soc.
\textbf{17}
(1977),
257--274.
%
\bibitem{Bb2012}
T. Bembom,
\textit{The capitulation problem in class field theory},
Dissertation,
Univ. G\"ottingen,
2012.
%
\bibitem{BEO2005}
H. U. Besche, B. Eick, and E. A. O'Brien,
\textit{The SmallGroups Library --- a Library of Groups of Small Order},
2005,
an accepted and refereed GAP package, available also in MAGMA.
%
\bibitem{BCP1997}
W. Bosma, J. Cannon, and C. Playoust,
\textit{The Magma algebra system. I. The user language}, 
J. Symbolic Comput.
\textbf{24}
(1997),
235--265.
%
\bibitem{BCFS2021}
W. Bosma, J. J. Cannon, C. Fieker, A. Steels (eds.),
\textit{Handbook of Magma functions},
Ed. 2.26,
Sydney,
2021.
%
\bibitem{BBH2017}
N. Boston, M. R. Bush and F. Hajir,
\textit{Heuristics for \(p\)-class towers of imaginary quadratic fields},
Math. Ann.
\textbf{368}
(2017),
No. 1,
633--669,
DOI 10.1007/s00208-016-1449-3.
%
\bibitem{BuMa2015}
M. R. Bush and D. C. Mayer,
\textit{\(3\)-class field towers of exact length \(3\)},
J. Number Theory
\textbf{147}
(2015),
766--777, \\
DOI 10.1016/j.jnt.2014.08.010.
%
\bibitem{ELNO2013}
B. Eick, C. R. Leedham-Green, M. F. Newman, and E. A. O'Brien,
\textit{On the classification of groups of prime-power order by coclass:
The \(3\)-groups of coclass \(2\)},
Int. J. Algebra Comput.
\textbf{23}
(2013),
1243--1288.
%
\bibitem{Fi2001}
C. Fieker,
\textit{Computing class fields via the Artin map},
Math. Comp.
\textbf{70}
(2001),
No. 235,
1293--1303.
%
\bibitem{GNO2006}
G. Gamble, W. Nickel, and E. A. O'Brien,
\textit{ANU \(p\)-Quotient --- \(p\)-Quotient and \(p\)-Group Generation Algorithms},
2006,
an accepted GAP package, available also in MAGMA.
%
\bibitem{HeSm1982}
F.-P. Heider und B. Schmithals,
\textit{Zur Kapitulation der Idealklassen in unverzweigten primzyklischen Erweiterungen},
J. Reine Angew. Math.
\textbf{336}
(1982),
1--25.
%
\bibitem{HEO2005}
D. F. Holt, B. Eick, and E. A. O'Brien,
\textit{Handbook of computational group theory},
Discrete mathematics and its applications,
Chapman and Hall/CRC Press,
Boca Raton,
2005.
%
\bibitem{KoVe1975}
H. Koch und B. B. Venkov,
\textit{\"Uber den \(p\)-Klassenk\"orperturm eines imagin\"ar-quadra\-tischen Zahlk\"orpers},
Ast\'erisque
\textbf{24--25}
(1975),
57--67.
%
\bibitem{MAGMA2021}
MAGMA Developer Group,
MAGMA \textit{Computational Algebra System},
Version 2.26-7,
Univ. Sydney,
2021, \\
\texttt{(http://magma.maths.usyd.edu.au)}.
%
\bibitem{Ma1991}
D. C. Mayer,
\textit{Principalization in complex \(S_3\)-fields},
Congressus Numerantium
\textbf{80}
(1991),
73--87
(Proc. of the Twentieth Manitoba Conf.
on Numerical Mathematics and Computing,
Winnipeg, Manitoba, Canada, 1990).
%
\bibitem{Ma2012a}
D. C. Mayer,
\textit{Transfers of metabelian \(p\)-groups},
Monatsh. Math.
\textbf{166}
(2012),
No. 3--4,
467--495, \\
DOI 10.1007/s00605-010-0277-x.
%
\bibitem{Ma2015a}
D. C. Mayer,
\textit{Periodic bifurcations in descendant trees of finite \(p\)-groups},
Adv. Pure Math.
\textbf{5}
(2015),
No. 1,
162--195,
DOI 10.4236/apm.2015.54020.
%
\bibitem{Ma2015b}
D. C. Mayer,
\textit{New number fields with known \(p\)-class tower}, 
Tatra Mt. Math. Pub.
\textbf{64}
(2015),
21--57, \\
DOI 10.1515/tmmp-2015-0040,
Special Issue on Number Theory and Cryptology \lq 15.
%
\bibitem{Ma2016a}
D. C. Mayer,
\textit{Artin transfer patterns on descendant trees of finite \(p\)-groups}, 
Adv. Pure Math.
\textbf{6}
(2016),
No. 2,
66--104,
DOI 10.4236/apm.2016.62008, 
Special Issue on Group Theory Research,
January 2016.
%
\bibitem{Ma2016b}
D. C. Mayer,
\textit{\(p\)-Capitulation over number fields with \(p\)-class rank two},
J. Appl. Math. Phys.
\textbf{4}
(2016),
No. 7,
1280--1293,
DOI 10.4236/jamp.2016.47135.
%
\bibitem{Ma2018}
D. C. Mayer,
\textit{Modeling rooted in-trees by finite \(p\)-groups},
Chapter 5, pp. 85--113,
in the Open Access Book \textit{Graph Theory --- Advanced Algorithms and Applications},
Ed. B. Sirmacek,
InTech d.o.o., Rijeka, January 2018,
DOI 10.5772/intechopen.68703.
%
\bibitem{Ma2020a}
D. C. Mayer,
\textit{Pattern recognition via Artin transfers applied to class field towers},
3rd International Conference on Mathematics and its Applications (ICMA) 2020,
Facult\'e des Sciences d' Ain Chock Casablanca (FSAC), Universit\'e Hassan II,
Casablanca, Morocco, invited keynote February 28, 2020, \\
\texttt{http://www.algebra.at/DCM@ICMA2020Casablanca.pdf}.
%
\bibitem{Ma2020b}
D. C. Mayer,
\textit{Schur \(\sigma\)-groups with abelian quotient invariants \((9,3)\)},
arXiv:2006.09177.
%
\bibitem{Ma2021}
D. C. Mayer,
\textit{Bicyclic commutator quotients with one non-elementary component},
arXiv:2108.10754.
%
\bibitem{Ne1989}
B. Nebelung,
\textit{Klassifikation metabelscher \(3\)-Gruppen
mit Faktorkommutatorgruppe vom Typ \((3,3)\)
und Anwendung auf das Kapitulationsproblem},
Inauguraldissertation,
Universit\"at zu K\"oln,
1989.
%
\bibitem{Nm1977}
M. F. Newman,
\textit{Determination of groups of prime-power order},
pp. 73--84
in: Group Theory, Canberra, 1975,
Lecture Notes in Math.,
Vol. \textbf{573}
(1977),
Springer,
Berlin.
%
\bibitem{Ob1990}
E. A. O'Brien, 
\textit{The p-group generation algorithm}, 
J. Symbolic Comput.
\textbf{9}
(1990),
677--698.
%
\bibitem{SoTa1934}
A. Scholz und O. Taussky,
\textit{Die Hauptideale der kubischen Klassenk\"orper imagin\"ar quadratischer Zahlk\"orper:
ihre rechnerische Bestimmung und ihr Einflu\ss\ auf den Klassenk\"orperturm},
J. Reine Angew. Math.
\textbf{171}
(1934),
19--41.
%
\bibitem{Sh1964}
I. R. Shafarevich,
\textit{Extensions with prescribed ramification points} (Russian),
Publ. Math., Inst. Hautes \'Etudes Sci. 
\textbf{18}
(1964),
71--95.
(English transl. by J. W. S. Cassels in
Amer. Math. Soc. Transl.,
II. Ser.,
\textbf{59}
(1966),
128--149.)
%
\bibitem{Ta1970}
O. Taussky,
\textit{A remark concerning Hilbert's Theorem \(94\)},
J. Reine Angew. Math.
\textbf{239/240}
(1970),
435--438.
%
\end{thebibliography}
\end{document}